\documentclass[aap,MSNbibl,dvips]{arximspdf}
\usepackage{mathrsfs}
\usepackage{graphicx}
\doi{10.1214/13-AAP981} %kopijuoti is PTS
\volume{24}
\issue{6}
\pubyear{2014}
\firstpage{2414}
\lastpage{2454}

\makeatletter
\newcommand{\widebar}{\overline}

\newcommand{\E}{\mathrm{E}}
\renewcommand{\P}{\mathrm{P}}
\renewcommand{\mid}{|}
\newcommand{\rrVert}{\Vert}
\newcommand{\rrvert}{\vert}
\newcommand{\llVert}{\Vert}
\newcommand{\llvert}{\vert}
\newcommand{\bi}{}
\newcommand{\binom}[2]{{{#1}\choose{#2}}}
\newcommand{\xrightarrow}[1]{\stackrel{#1}{\longrightarrow}}
\newcommand{\argmin}{\mathop{\arg\min}}
\newcommand{\eqdist}{\stackrel{\mathrm{D}}{=}} %equality in
\newtheorem{theorem}{Theorem}
\newtheorem{lemma}{Lemma}
\newproclaim{defnn}{Definition}
\makeatother

\begin{document}
\begin{frontmatter}

\title{Belief propagation for optimal edge cover in the random
complete graph}
\runtitle{Belief propagation for edge cover}

\begin{aug}
\author{\fnms{Mustafa} \snm{Khandwawala}\corref{}\ead[label=e1]{mustafa@ece.iisc.ernet.in}\thanksref{t1}}
\and
\author{\fnms{Rajesh} \snm{Sundaresan}\ead[label=e2]{rajeshs@ece.iisc.ernet.in}\thanksref{t1,t2}}
\runauthor{M. Khandwawala and R. Sundaresan}
\affiliation{Indian Institute of Science}
\address{Department of Electrical Communication Engineering\\
Indian Institute of Science\\
Bangalore 560012\\
India\\
\printead{e1}\\
\phantom{E-mail:\ }\printead*{e2}} %adresu isvedimo komanda gale!
\end{aug}
\thankstext{t1}{Supported by the Department of Science and Technology,
Government of India, by a TCS fellowship grant.}
\thankstext{t2}{Supported by an Indo-US Science and Technology Forum research fellowship grant.}

% HISTORY:
\received{\smonth{12} \syear{2012}}
\revised{\smonth{10} \syear{2013}}

% ABSTRACT
%
\begin{abstract}
We apply the objective method of Aldous to the problem of finding the
minimum-cost edge cover of the complete graph with random independent
and identically distributed edge costs. The limit, as the number of
vertices goes to infinity, of the expected minimum cost for this
problem is known via a combinatorial approach of Hessler and W\"
{a}stlund. We provide a proof of this result using the machinery of the
objective method and local weak convergence, which was used to prove
the $\zeta(2)$ limit of the random assignment problem. A proof via the
objective method is useful because it provides us with more information
on the nature of the edge's incident on a typical root in the minimum-cost edge cover. We further show that a belief propagation algorithm
converges asymptotically to the optimal solution. This can be applied
in a computational linguistics problem of semantic projection. The
belief propagation algorithm yields a near optimal solution with lesser
complexity than the known best algorithms designed for optimality in
worst-case settings.
\end{abstract}

% KEYWORDS
\begin{keyword}[class=AMS]
\kwd[Primary ]{60C05}
\kwd[; secondary ]{68Q87}
\kwd{82B44}
\end{keyword}
\begin{keyword}
\kwd{Belief propagation}
\kwd{edge cover}
\kwd{local weak convergence}
\kwd{objective method}
\kwd{semantic projection}
\end{keyword}

\end{frontmatter}

%s1 #&#
%s1 ###
\section{Introduction} \label{secintro}
Suppose that we are given a graph $G$ with vertex set $V$ and edge set
$E$, denoted $G = (V,E)$. Each edge $e \in E$ has a weight $\xi_e \in
\mathbf{R}_{+}$. Alternatively, we are given a bipartite graph with a vertex
set $V = V_1 \cup V_2$, a union of two disjoint vertex subsets, and an
edge set $E \subset V_1 \times V_2$. An \textit{edge cover} for the graph
is a subset of edges that hits (covers) every vertex. The cost of an
edge cover is the sum of the weights of edges in the cover. Our
interest in this paper is on minimum-cost edge covers on the complete
graph (denoted $K_n$ when $|V| = n$) and on the complete bipartite
graph (denoted $K_{n,n}$ when $|V_1| = |V_2| = n$), when the edge
weights are independent random variables, each with the exponential
distribution of mean 1.

The following example on a bipartite graph illustrates how minimum-cost
edge covers arise in practice.

\subsection*{An example of semantic projection}
Computational linguists have recently been interested in machine-based
natural language processing. These include part-of-speech tagging,
parsing, and at a higher level, semantic role parsing \cite
{PadLap2006} which, for example, would enable an automatic recognition
that the sentences ``Mary sold the book to John'' and ``The book was
sold by Mary to John'' have the same semantic roles. (This example is
taken from Wikipedia \cite{WikiSemanticrolelabeling}.) Currently,
English is blessed with the availability of a large amount of annotated
texts as training data while most others languages lack this advantage.
Semantic projection exploits the availability of (1) parallel corpora
of translated texts and (2) higher quality parsing tools in one
language in order to transfer annotations from the resource-rich
language to the other.

Pad\'{o} and Lapata \cite{PadLap2006} provide one method to do this
where a minimum-cost edge cover naturally arises. The source and target
sentences in the two languages are first broken into \textit{linguistic
units} to yield sets $V_1$ and $V_2$ of the respective linguistic
units. These linguistic units are then viewed as vertices of a complete
bipartite graph. Let $R$ be some finite set of \textit{semantic roles},
which can be viewed for our purposes as abstract annotations. The
parsing tool on the source side is used to find a semantic role
assignment $\textsf{role}_1\dvtx  R \rightarrow2^{V_1}$, where the
subscript refers to the source language. A \textit{dissimilarity} measure
based on linguistic considerations is then assigned to every pair of
linguistic units across the languages and is denoted $\xi\dvtx  V_1 \times
V_2 \rightarrow\mathbf{R}_{+}$. A~decision procedure uses these dissimilarity
scores to find a subset $\mathcal{C} \subset V_1 \times V_2$ of
semantically aligned units. Pad\'{o} and Lapata \cite{PadLap2006}
argue that a minimum-cost edge cover is a good choice for this semantic
alignment. It allows a linguistic unit in one language (an element of
say $V_2$) to map to several units in the other language (a subset of
$V_1$), and vice-versa. For example, the linguistic units ``to be on
time'' and ``punctual'' (English) could both be mapped with small, but
possibly different, dissimilarity scores to ``p\"{u}nktlich'' (German),
and both edges may be picked by a good candidate edge cover. The
covering property of the edge cover enables all source and target
vertices to participate and thus has the potential to capture important
connections between linguistic units, which may otherwise be missed.
The minimum cost property attempts to provide an economical semantic
alignment and further captures global alignments as compared to
previously proposed local decision procedures. Once the minimum-cost
edge cover is found by the decision procedure, semantic roles are then
assigned on the target side as
\[
\textsf{role}_2(r) = \bigl\{ j \mid\mbox{ there is an } i \in\textsf
{role}_1(r) \mbox{ such that } (i,j) \in\mathcal{C} \bigr\}.
\]
Pad\'{o} and Lapata \cite{PadLap2006} compare the goodness of their
decision procedures based on minimum-cost edge cover (and perfect
matching) with some other prior approaches on a data set of about 1000
sentences. Real data sets are of course much larger. The resulting
graph, when restricted to edges of small weight (i.e., edges signifying
low dissimilarity and therefore good correspondence), can be modeled as
a large, but sparse, random graph. If $|V_1| = O(|V_2|) = n$,
algorithms used by Pad\'{o} and Lapata \cite{PadLap2006} to find the
minimum-cost edge cover take $O(n^3)$ operations, in the worst case.

The actual results of the Pad\'{o} and Lapata \cite{PadLap2006}
experiments need not concern us here. For a list of challenges that
arise in the implementation of the above approach and methods to address
them, we refer the linguistically inclined reader to~\cite{PadLap2009} and
references therein. What we shall take with us as we move forward are
the observations that (1) edge covers arise in practice on large graphs
that can be modeled by sparse random graphs and (2) algorithmic
simplifications that reduce complexity are of practical value.

We shall for simplicity focus on minimum-cost edge covers on the
complete graph $K_n$ on $n$ vertices. All our results carry over to
$K_{n,n}$ with only scaling factor modifications.
Recall that the edge capacities are independent, each edge having the
exponential distribution with mean 1. This is a
typical \textit{mean-field} model which captures sparsity of the graph
depicting linguistic units and associated edges in the above example, but
ignores correlations among edge weights. See Section~\ref{secsummary}
for another
geometric setting where the same mean field models arise. Let $C_n$ be
the cost of the minimum-cost edge cover of $K_n$. We prove that the
expected value of $C_n$ converges to the constant $W(1)+W(1)^2/2$,
which is approximately 0.728. (The function $W(\cdot)$ is Lambert's
$W$-function, which is the inverse of $f\dvtx  [0,\infty) \rightarrow
[0,\infty), f(x) = xe^x$; $W(1) \approx0.567$.) Further, and more
importantly from an application perspective, we show that a belief\vspace*{1pt}
propagation algorithm can be used to find asymptotically optimal edge
covers in $O(n^2)$ steps. The results, with only scaling factor
changes, hold for the complete bipartite graphs $K_{n,n}$.

The result regarding the limit on $K_{n,n}$ has been proved before by
Hessler and W\"{a}stlund in \cite{HesWas2010} using a combinatorial
approach. A proof based on a game formulation is contained in \cite
{Was2009}. We discuss these works at the end of this \hyperref[secintro]{Introduction}. Our
focus in this article is on using the \textit{objective method} for this
problem and on devising a belief propagation algorithm.

The roots of the objective method lie in Aldous's 1992 paper \cite
{Ald1992} on the assignment problem. The problem of finding the minimum
cost matching on the complete bipartite graph with independent and
identically distributed edge costs, termed as the random assignment
problem in literature, inspired a series of works in combinatorial
probability. M\'{e}zard and Parisi \cite{MezPar1987}, using the \textit{cavity method} of statistical physics, conjectured in 1987 that the
expected minimum cost for the random assignment problem on the
bipartite graph $K_{n,n}$ converges to $\zeta(2)=\sum_{k=1}^{\infty}
k^{-2}$ as
$n$ goes to infinity. This was proved rigorously by Aldous \cite
{Ald2001} in 2001 by extending the proof of existence of the limit
contained in \cite{Ald1992}. Several other proofs have been provided
for the limit in subsequent works. %need references
%here***********************

\if0
The key ingredient of the proof in \cite{Ald2001} is to obtain a
distributional identity to construct a process which serves as a
suitable limit for the optimization problem. This allows one to relate
quantities in the finite graph with appropriately-defined quantities on
the infinite limit object. For instance, the following quantities are
shown to converge to corresponding quantities on a limit infinite graph
in \cite{Ald2001}
%
\begin{itemize}
\fi

In \cite{Ald2001}, Aldous related the problem on $K_{n,n}$ to one on a
suitable limit object. Several calculations become easier on the limit
object. In this case, the limit is a tree, the so-called Poisson
weighted infinite tree or PWIT, with many useful symmetries. Aldous
used these symmetries to construct a \textit{distributional identity},
that then served as a guide for solving the random assignment problem
rigorously. With this approach, Aldous showed that the following
quantities converge to the corresponding quantities on the limit object:
\begin{itemize}
\item the expected cost of optimal matching on $K_{n,n}$;
\item the distribution of the cost of the matching edge incident on a
typical node of~$K_{n,n}$;
\item the probability that the matching edge incident on a typical node
of $K_{n,n}$ is the $k$th smallest of all the edges incident on it.
\end{itemize}
It turns out that the limit object, and hence the answers, remain the
same for problems on the complete bipartite graph $K_{n,n}$ and on the
complete graph $K_n$. One dividend of a proof via the objective method
is that we have answers to several ancillary questions such as the
second and third bullets above. The ability of the objective method to
provide these auxiliary results motivates us to solve the problem of
optimal edge cover via the objective method.\looseness=-1

From an algorithms perspective, the cavity equations suggest a natural
iterative decentralized message passing algorithm, some versions of
which are commonly called belief propagation (BP) in the computer
science literature. For many combinatorial optimization problems, a BP
algorithm can be set up to converge to the correct solution on graphs
without cycles. Bayati, Shah and Sharma \cite{BayShaSha2008} proved
that the BP algorithm for maximum weight matching on bipartite graphs
converges to the correct value as long as the maximum weight matching
is unique. Salez and Shah \cite{SalSha2009} studied the random
assignment problem and proved a tighter connection with the limit
object. They showed that that a BP algorithm on $K_{n,n}$ converges to an
update rule on the limit PWIT of \cite{Ald2001}. The iterates on the
limit graph converge in distribution to the minimum cost assignment.
The iterates are near the optimal solution in $O(n^2)$ steps, whereas
the worst case optimal algorithm on bipartite graphs is $O(n^3)$
[expected time $O(n^2 \log n)$ for i.i.d. edge capacities]; see Salez
and Shah \cite{SalSha2009} and references therein. We show a similar
complexity improvement for the edge-cover problem.

The objective method is quite powerful to be applicable to several
combinatorial probability problems. See Aldous and Steele \cite
{AldSte2004} for a survey. Aldous and Bandopadhyay \cite{AldBan2005}, Section~7.5, outline the steps of Aldous's program to
establish the validity of the cavity method, which we quote in
Section~\ref{secsummary}. However, each problem requires specific
proofs, and we are still far from a complete theory applicable to a
wide class of problems. The edge-cover problem itself poses some modest
problem-specific challenges which we overcome in this paper. These
include (1) a proof of existence and uniqueness of a solution to the
distributional identity associated with the edge-cover problem, (2) a
proof of a property called \textit{endogeny} of a process on the tree
associated with the distributional identity, (3)~a~proof of optimality
of the edge-cover selection on the PWIT as suggested by the
distributional identity and eventually (4) a proof that a BP algorithm
converges to an asymptotically optimal edge cover on the random
complete graph. See Section~\ref{secsummary} for a more detailed summary.

Before we end this introduction, we would like to mention two other
approaches that have been used to solve related combinatorial
optimization problems, in particular, matching, edge cover and
travelling salesman problems. One approach used by W\"{a}stlund in
\cite{Was2009,Was2012} calls for a ``boundary conditioning'' parameter
to study ``diluted'' versions of the optimization problems, eventually
driving the parameter to infinity, and thereby relating the resulting
limiting problem with the undiluted versions. For example, in the
matching case, diluted matching is a partial matching with each
unmatched vertex paying a cost equal to the parameter. W\"{a}stlund
then formulates the optimization problem in terms of a game played on
the graph. A second and more combinatorial approach is used by W\"
{a}stlund in \cite{Was2010} for matching and TSP and in \cite
{HesWas2010} for the edge-cover problem. These works study the
respective optimization problems as certain flow problems on bipartite
graphs. The feasible solutions to these flow problems have a fixed
number of edges $k$. A~recursive relation on $k$ is obtained for the
cost of the optimal solution. As our focus is on the objective method,
we do not dwell any more on these approaches.

%s2 #&#
%s2 ###
\section{Main results} \label{secresults}
Our first result establishes the limit of the expected minimum cost of
the random edge-cover problem.
%
%th1 #&#
\begin{theorem} \label{thmlimit}
On $K_n$, we have
%
%e1 #&#
%e1 ###
\begin{equation}
\lim_{n\rightarrow\infty}\E{C_n}=W(1)+\frac{W(1)^2}{2}.
\end{equation}
\end{theorem}

Our second result shows that a belief propagation algorithm gives an
edge cover that is asymptotically optimal as $n\rightarrow\infty$. We will
use the result that the update rule of BP converges to an update rule
on a limit infinite tree. For this we define\vadjust{\goodbreak} the BP algorithm on an
arbitrary graph $G=(V,E)$ with edge costs. For an edge $e=\{v,w\}\in
E$, we write its cost as $\xi_{G} (e )$ or $\xi_{G}
(v,w )$. For each
vertex $v\in V$, we associate a nonempty subset of its neighbors $\pi
_{G}^{}(v)$. By taking a union of all edges of the form $\{v,w\}, w\in
\pi_{G}^{}(v)$, we get an edge cover of $G$ which we will denote by
$\mathcal{C}(\pi_{G}^{})$.

The BP algorithm is an iterative message passing algorithm. In each
iteration $k\ge0$, every vertex $v\in V$ sends a message
$X_{G}^{k} (w,v )$ to each neighbor $w\sim v$ according to
the following rules:
\begin{description}
\item[Initialization:]%
%e2 #&#
%e2 ###
\begin{equation}
\label{eqbpInit} X_{G}^{0} (w,v )=0.
\end{equation}
\item[Update rule:]%
%e3 #&#
%e3 ###
\begin{equation}
\label{eqbpUpdate} X_{G}^{k+1} (w,v )=\min_{{u\sim v,u\neq w}}
\bigl\{ \bigl(\xi_{G} (v,u )-X_{G}^{k} (v,u )
\bigr)^+ \bigr\}.
\end{equation}
\item[Decision rule:]%
%e4 #&#
%e5 #&#
%e5 ###
%e4 ###
\begin{eqnarray}
\label{eqbpDec} \pi_{G}^{k}(v)&=&\argmin_{u\sim v} \bigl\{
\bigl(\xi_{G} (v,u )-X_{G}^{k} (v,u ) \bigr)^+
\bigr\},
\\
\mbox{Edge cover}&=&\mathcal{C}\bigl(\pi_{G}^{k}(v)\bigr).
\end{eqnarray}
\end{description}

We analyze the belief propagation algorithm for $G=K_n$ and i.i.d.
exponential random edge costs, and prove that after a sufficiently
large number of iterates, the expected cost of the assignment given by
the BP algorithm is close to the limit value in Theorem~\ref{thmlimit}.
%
%th2 #&#
\begin{theorem} \label{thmbpconv}
On $K_n$, we have
%
%e6 #&#
%e6 ###
\begin{equation}
\label{eqthmbpconv} \lim_{k\rightarrow\infty}\lim_{n\rightarrow\infty}\E{ \biggl[
\sum_{e\in\mathcal{C}(\pi
_{K_n}^{k})}\xi_{K_n} (e ) \biggr]}=W(1)+
\frac{W(1)^2}{2}.
\end{equation}
\end{theorem}

The formal statements on the bipartite complete graph $K_{n,n}$ with
i.i.d. expontial distribution of mean 1 are the following and are
stated without proof.\looseness=-1

%th3 #&#
\begin{theorem} On $K_{n,n}$, we have
%
%e7 #&#
%e7 ###
\begin{equation}
\lim_{n\rightarrow\infty}\E{C_n}=2W(1)+W(1)^2.
\end{equation}
\end{theorem}
%
%
%th4 #&#
\begin{theorem} On $K_{n,n}$, we have
%
%e8 #&#
%e8 ###
\begin{equation}
\lim_{k\rightarrow\infty}\lim_{n\rightarrow\infty}\E{ \biggl[\sum
_{e\in\mathcal{C}(\pi_{K_{n,n}}^{k})}\xi_{K_{n,n}} (e ) \biggr]}=2W(1)+W(1)^2.
\end{equation}
\end{theorem}

%s3 #&#
%s3 ###
\section{Local weak convergence} \label{seclwConv}
In this section, we recollect the terminology for defining convergence
of graphs.

%s3.1 #&#
%s3.1 ###
\subsection{Rooted geometric networks}
A graph $G=(V,E)$ along with a length function $l\dvtx E\rightarrow
(0,\infty]$ is
called a \textit{network}.\vadjust{\goodbreak} The \textit{distance} between two vertices in
the
network is the infimum of the sum of lengths of the edges of a path
connecting the two vertices, the infimum being taken over all such
paths. We call the network a \textit{geometric network} if for each vertex
$v\in V$ and positive real $\rho$, the number of vertices within a
distance $\rho$ of $v$ is finite. We denote the space of geometric
networks by ${\mathcal{G}}$.

A geometric network with a distinguished vertex $v$ is called a \textit{rooted geometric network} with root $v$. We denote the space of all
connected rooted geometric networks by ${\mathcal{G}_*}$. In
${\mathcal{G}_*}$ we do
not distinguish between rooted isomorphisms of the same network. We
will use the notation $(G,o)$ to denote an element of ${\mathcal
{G}_*}$ which
is the isomorphism class of rooted networks with underlying network $G$
and root~$o$.

%s3.2 #&#
%s3.2 ###
\subsection{Local weak convergence}
We call a positive real number $\rho$ a \textit{continuity point} of $G$
if no vertex of $G$ is exactly at a distance $\rho$ from the root of
$G$. Let $\mathcal{N}_\rho(G)$ denote the neighborhood of the root of
$G$ up to
distance $\rho$. $\mathcal{N}_\rho(G)$ contains all vertices of $G$
which are
within a distance $\rho$ from the root of $G$ (Figure~\ref{fignbd}).
We take $\mathcal{N}_\rho(G)$ to be an element of ${\mathcal{G}_*}$
by inheriting the
same length function $l$ as $G$, and the same root as that of $G$.

%
%f1 #&#
%f1 ###
\begin{figure}%[b]

\includegraphics{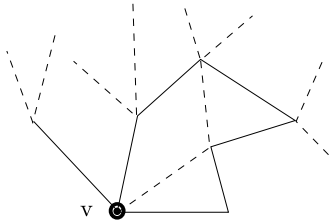}

\caption{Neighborhood $\mathcal{N}_\rho(G)$ of graph $G$. The solid
edges form
the neighborhood, and form paths of length at most $\rho$ from the
root $v$. Dashed edges are the other edges of $G$.}\label{fignbd}
\end{figure}

We say that a sequence of rooted geometric networks $G_n, n\ge1$, \textit{converges locally} to an element $G_\infty$ in ${\mathcal{G}_*}$ if
for each
continuity point $\rho$ of $G_\infty$, there is an $n_\rho$ such
that for all $n\ge n_\rho$, there exists a graph isomorphism $\gamma
_{n,\rho}$ from $\mathcal{N}_\rho(G_\infty)$ to $\mathcal{N}_\rho
(G_n)$ that maps the root
of the former to the root of the latter, and for each edge~$e$ of
$\mathcal{N}_\rho(G_\infty)$, the length of $\gamma_{n,\rho}(e)$
converges to the
length of $e$ as $n\rightarrow\infty$.

\if0
It is possible to define weak convergence of measures on ${\mathcal{G}_*}$,
and term it \textit{local weak convergence}.
\fi

The space ${\mathcal{G}_*}$ can be suitably metrized to make it a separable
and complete metric space. One can then consider probability measures
on this space and endow that space with the topology of weak
convergence of measures. This notion of convergence is called \textit{local weak convergence}.

In our setting of complete graphs $K_n=(V_n,E_n)$ with random i.i.d.
edge costs $ \{\xi_e,e\in E_n \}$, we regard the edge costs
to be the
lengths of the edges, and declare a vertex of $K_n$ chosen uniformly
at random as the root of $K_n$. This makes $K_n$ along with its root
a random element of ${\mathcal{G}_*}$. We rescale the edge costs such
that for
each~$n$, $ \{\xi_e,e\in E_n \}$ are i.i.d. random
variables with
mean $n$ exponential distribution. We will denote this random, rooted,
rescaled version of the $n$-vertex complete graph by $\widebar{K}_n$ to
distinguish it from the $K_n$ defined earlier. Theorem~\ref
{thmpwitConv} stated below (from \cite{Ald1992}) says
that the sequence of random geometric networks $\widebar{K}_n$
converges in the
local weak sense to an element of ${\mathcal{G}_*}$ called the \textit{Poisson
weighted infinite tree} (\textit{PWIT}).

%s3.3 #&#
%s3.3 ###
\subsection{Poisson weighted infinite tree}
We use the notation from \cite{SalSha2009} to define the~PWIT.

Denote by $\mathcal{V}$ the set of all finite words over the alphabet
$\mathbf{N}= \{1,2,3,\ldots \}$. Let $\phi$ denote the
empty string
and ``.'' the concatenation operator. For any $v\in\mathcal{V}$ write
$\llvert v\rrvert $ for the length of string $v$, and if $v\neq\phi$
write $\dot
{v}$ for the string obtained by removing the last letter of $v$.

Construct an undirected graph $\mathcal{T}=(\mathcal{V},\mathcal{E})$ on
$\mathcal{V}$ with
the edge set
\[
\mathcal{E}= \bigl\{ \{v,v.i \},v\in\mathcal{V},i\in \mathbf{N} \bigr\}.
\]
Set $\phi$ to be the root of $\mathcal{T}$. Then $\mathcal{T}$ is an infinite
rooted tree with each vertex having a countably infinite number of
children. Construct a family of independent Poisson processes of
intensity 1 on $\mathbf{R}_{+}\dvtx  \{\bolds{\xi}^v=(\xi
^v_1,\xi^v_2,\ldots),v\in\mathcal{V} \}$. Assign to each
edge $ \{v,v.i \}$ in $\mathcal{E}$ the
length $\xi^v_i$. $\mathcal{T}$ is then a random element of
${\mathcal{G}_*}$, and
we call it the Poisson weighted infinite tree (PWIT) (Figure~\ref{figpwit}).

%
%f2 #&#
%f2 ###
\begin{figure}[b]

\includegraphics{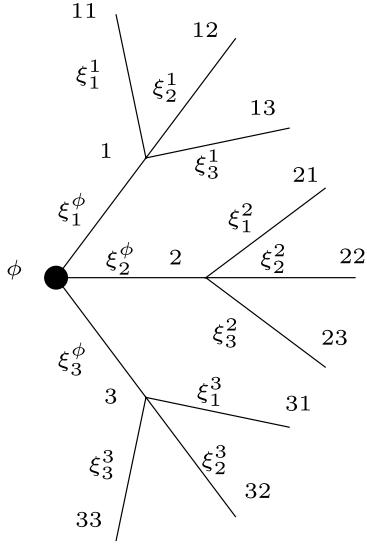}

\caption{PWIT $\mathcal{T}$ up to depth 2, with only the first three
children of each vertex shown.}\label{figpwit}
\end{figure}

%
%th5 #&#
\begin{theorem}[(\cite{Ald1992})] \label{thmpwitConv}
The sequence of uniformly rooted random networks $\widebar{K}_n$
converges to
the PWIT $\mathcal{T}$ as $n\rightarrow\infty$ in the sense of local weak
convergence.
\end{theorem}
A similar result was earlier established by Hajek \cite{Haj1990}, Section~IV, for a class of sparse Erd\H{o}s--R\'{e}nyi
random graphs.
The above theorem says that if we look at an arbitrary, large but fixed
neighborhood of the root of $\widebar{K}_n$, then for large $n$ it
looks like
the corresponding neighborhood of the root of $\mathcal{T}$. This suggests
that if boundary conditions can be ignored, we may be able to relate
optimal edge covers on $\widebar{K}_n$ with an appropriate edge cover
on $\mathcal{T}
$ [to be precise, an optimal \textit{involution invariant} edge cover
(Section~\ref{secinvInv}) on the PWIT]. Furthermore, the local
neighborhood of the root of $\widebar{K}_n$ is a tree for large enough $n$
(with high probability). So we may expect belief propagation on
$\widebar{K}_n$
to converge. Both the above observations are true in the matching case;
the former was established in \cite{Ald1992,Ald2001}, and the latter
was shown in~\cite{SalSha2009}. We now extend these ideas to prove
similar results for the edge-cover problem.
%s4 #&#
%s4 ###
\section{Recursive distributional equation} \label{secrde}
%s4.1 #&#
%s4.1 ###
\subsection{A heuristic recursion}
The PWIT $\mathcal{T}$ is an infinite graph, and it is clear that any edge
cover on it must have infinite cost. So it does not make sense to talk
about a minimum-cost edge cover on $\mathcal{T}$. However, for a
moment let
us pretend to perform operations on the minimum cost as if it were a
finite quantity. Write $C(\mathcal{T})$ for this minimum cost, and define
%
%e9 #&#
%e9 ###
\begin{equation}
\label{eqdiff} D(\mathcal{T})=\bigl(C(\mathcal{T})-C\bigl(\mathcal{T}\setminus
\{\phi \}\bigr)\bigr)^+,
\end{equation}
where $C(\mathcal{T}\setminus \{\phi \})$ is the minimum
cost of edge
cover on the subgraph of $\mathcal{T}$ obtained by removing the root. Note
that $D(\mathcal{T})$ denotes the difference between the minimum cost of
edge cover of $\mathcal{T}$ and the minimum cost of partial edge cover of~$\mathcal{T}$ where the root $\phi$ can be left uncovered.

%
%f3 #&#
%f3 ###
\begin{figure}

\includegraphics{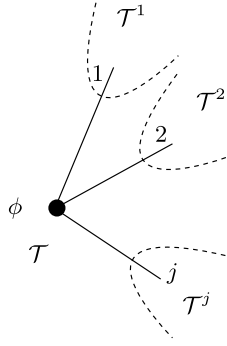}

\caption{PWIT $\mathcal{T}$ with the subtrees $\mathcal{T}^j$ at node
$j$.}\label{figsubtrees}
\end{figure}

If $j$ is a child of the root, let $\mathcal{T}^j$ denote the induced
subgraph of $\mathcal{T}$ containing~$j$ and all its descendants, and view
it as a rooted network with root $j$ (Figure~\ref{figsubtrees}).
Define $D(\mathcal{T}^j)$ accordingly, and observe from the symmetry of
$\mathcal{T}$ that $ \{D(\mathcal{T}^j),j\ge1 \}$ are
i.i.d., and have the same
distribution as $D(\mathcal{T})$. We give a heuristic argument that
$D(\mathcal{T}
)$ satisfies the following relation:
%
%e10 #&#
%e10 ###
\begin{equation}
\label{eqDT} D(\mathcal{T})=\min_{j\ge1}\bigl(
\xi^\phi_j-D\bigl(\mathcal{T}^j\bigr)\bigr)^+.
\end{equation}

We can write $C(\mathcal{T}\setminus \{\phi \})$ in terms
of edge covers
on the subtrees $\mathcal{T}^j,j\ge1$, as
%
%e11 #&#
%e11 ###
\begin{equation}
C\bigl(\mathcal{T}\setminus \{\phi \}\bigr)=\sum_{j\in\mathbf
{N}}C
\bigl(\mathcal{T}^j\bigr).
\end{equation}
Let us consider edge covers in which the edges covering the root are
incident on the vertices in a fixed subset $A$ of the children of the
root. The minimum cost among such edge covers can be written as
\[
\sum_{j\in A} \bigl(\xi^\phi_j+
\min \bigl\{C\bigl(\mathcal {T}^j\bigr),C\bigl(\mathcal{T}
^j\setminus \{j \}\bigr) \bigr\} \bigr)+\sum
_{i\in\mathbf{N}\setminus A}C\bigl(\mathcal{T}^i\bigr).
\]
$C(\mathcal{T})$ is the minimum of the above value taken over all nonempty
$A$, that is,
%
%e12 #&#
%e12 ###
\begin{equation}
C(\mathcal{T}) = \min_{A\ \mathrm{nonempty}} \biggl\{\sum
_{j\in
A} \bigl(\xi^\phi_j+\min \bigl\{C
\bigl(\mathcal{T}^j\bigr),C\bigl(\mathcal {T}^j\setminus
\{j \}\bigr) \bigr\} \bigr)+\sum_{i\in\mathbf
{N}\setminus A}C\bigl(
\mathcal{T}^i\bigr) \biggr\}.\hspace*{-35pt}
\end{equation}
Thus we can write
\begin{eqnarray*}
D(\mathcal{T})&=& \biggl(\min_{\llvert A\rrvert \ge1}\sum
_{j\in
A} \bigl(\xi^\phi _j-\bigl(C\bigl(
\mathcal{T}^j\bigr)-C\bigl(\mathcal{T}^j\setminus \{j \}
\bigr)\bigr)^+ \bigr) \biggr)^+
\\
&=& \biggl(\min_{\llvert A\rrvert \ge1}\sum_{j\in A}
\bigl(\xi^\phi _j-D\bigl(\mathcal{T}^j\bigr)
\bigr) \biggr)^+.
\end{eqnarray*}
To minimize the term within parentheses, we must include all those
indices $j$ for which the summand $(\xi_j^{\phi} - D(\mathcal
{T}^j))$ is
negative. If the terms are positive for all indices $j$, $A$ must be
the singleton where the minimum is attained among all indices. By then
taking the positive part, equation~(\ref{eqDT}) follows.

Although $D(\mathcal{T})$ and $D(\mathcal{T}^j)$ are not well-defined
quantities,
we shall prove that there is a nonnegative random variable $X$ and
i.i.d. random variables $X_j,j\ge1$, having the same distribution as
$X$, such that
%
%e13 #&#
%e13 ###
\begin{equation}
\label{eqECrde} X=\min_{j\ge1}(\xi_j-X_j)^+,
\end{equation}
where $ \{\xi_j,j\ge1 \}$ are points of a Poisson process
of rate 1
on $\mathbf{R}_{+}$, independent of $ \{X_j,j\ge1 \}$.

%s4.2 #&#
%s4.2 ###
\subsection{Recursive distributional equations and recursive tree processes}
Equations of the form of (\ref{eqECrde}) are termed as recursive
distributional equations in \cite{AldBan2005}. Specifically, if
$\mathcal{P}(S)$ denotes the space of probability measures on a space
$S$, a~recursive distributional equation (RDE) is a fixed-point
equation on $\mathcal{P}(S)$ of the form
%
%e14 #&#
%e14 ###
\begin{equation}
\label{eqrde} X\eqdist g\bigl(\xi;(X_j,1\le j< N)\bigr),
\end{equation}
where $X_j,j\ge1$ are i.i.d. $S$-valued random variables having the
same distribution as $X$, and are independent of the pair $(\xi,N)$,
$\xi$ is a random variable on some space and $N$ is a random variable
on $\mathbf{N}\cup \{+\infty \}$. $g$ is a given
$S$-valued function.
A~solution to the RDE is a common distribution of $X, X_j,j\ge1$,
satisfying (\ref{eqrde}).

We can use relation (\ref{eqrde}) to construct a tree indexed
stochastic process, say $X_{\bi i},{\bi i}\in\mathcal{V}$, which is called
a recursive tree process (RTP) \cite{AldBan2005}. Associate to each
vertex ${\bi i}\in\mathcal{V}$, an independent copy $(\xi_{\bi
i},N_{\bi
i})$ of the pair $(\xi,N)$, and require $X_{\bi i}$ to satisfy
\[
X_{\bi i}\eqdist g\bigl(\xi_{\bi i};(X_{{\bi i}.j},1\le j<
N_{\bi i})\bigr)
\]
with $X_{\bi i}$ independent of $ \{(\xi_{\bi i^\prime},N_{\bi
i^\prime}){ \vert}\llvert \bi i^\prime\rrvert <\llvert \bi i\rrvert  \}$.
If $\mu\in\mathcal{P}(S)$ is a solution to the RDE~(\ref{eqrde}),
there exists a stationary RTP; that is, each $X_{\bi i}$ is distributed
as $\mu$. Such a process is called an invariant RTP with marginal
distribution $\mu$.

%s4.3 #&#
%s4.3 ###
\subsection{Solution to the edge cover RDE}%
%th6 #&#
\begin{theorem}\label{thmRDEsolution}
The unique solution to the RDE~(\ref{eqrde}) is the c.d.f. $F_*$ whose
complementary c.d.f. $\widebar{F}_*$ is given by
%
%e15 #&#
%e15 ###
\begin{equation}
\label{eqfpdist} \widebar{F}_*(y)=\cases{ W(1)e^{-y}, &\quad if $y\ge0$,
\cr
1,&\quad if $y<0$.}
\end{equation}
The function $W$ above is Lambert's $W$-function, the inverse of
$f\dvtx \mathbf{R}_{+}\rightarrow\mathbf{R}_{+},f(x)=xe^x$. In
particular, $W(1)e^{W(1)}=1$.
\end{theorem}
\begin{pf}
Let $\mu$ be a solution to the RDE (\ref{eqECrde}), and let $F$ be
its c.d.f. Take $X_j,j\ge1$ i.i.d. with distribution $\mu$. Then
$ \{(\xi_j,X_j),j\ge1 \}$ is a Poisson process on $\mathbf
{R}_{+}\times\mathbf{R}_{+}$
with intensity $\mathrm{d}z \, \mathrm{d}F(x)$.
For $y\in\mathbf{R}_{+}$,
\begin{eqnarray*}
\P{(X>y)}&=&\P{ \Bigl(\min_{j\ge1}(\xi_j-X_j)^+>y
\Bigr)}
\\
&=&\P{ \bigl(\mbox{No point of } \bigl\{(\xi_j,X_j) \bigr
\}\mbox{ in } \bigl\{(z,x)\dvtx z-x\le y \bigr\} \bigr)}
\\
&=&\exp \biggl(-\int_{z=0}^y\int
_{x=0}^\infty \mathrm{d}F(x)\,
\mathrm{d}z-\int_{z=y}^\infty\int_{x=z-y}^\infty
\mathrm{d}F(x)\,\mathrm{d}z \biggr)
\\
&=&e^{-y}\exp \biggl(-\int_{t=0}^\infty\int
_{x=t}^\infty \mathrm{d}F(x)\,
\mathrm{d} t \biggr)
\\
&=&e^{-y}\exp \biggl(-\int_{0}^\infty
\bigl(1-F(t)\bigr)\, \mathrm{d}t \biggr).
\end{eqnarray*}
Writing $\widebar{F}(t)=1-F(t)$, we have
\[
\widebar{F}(y)=e^{-y}\exp \biggl(-\int_{0}^\infty
\widebar{F}(t)\, \mathrm{d}t \biggr)\qquad\mbox{for all }y\ge0.
\]
Let $c=\exp (-\int_{0}^\infty\widebar{F}(t)\, \mathrm{d}t )$. Then,
using $\widebar{F}(t)=ce^{-t}$ in the expression for $c$ gives
\[
c=\exp \biggl(-\int_{0}^\infty ce^{-t}
\mathrm{d}t \biggr)=e^{-c}.
\]
The unique $c$ satisfying the above equation is $c=W(1)$. This proves
that $F$ must be the c.d.f. $F_*$.
\end{pf}
%
%s5 #&#
%s5 ###
\section{Unimodularity and involution invariance} \label{secinvInv}
In Section~\ref{seclwConv} we defined the space ${\mathcal{G}_*}$ as
the set
of connected rooted geometric networks. Now define ${\mathcal
{G}_{**}}$ as the
space of connected geometric networks with an ordered pair of
distinguished vertices. Again, we do not distinguish between
isomorphisms in ${\mathcal{G}_{**}}$, and denote by $(G,o,x)$ the isomorphism
class of elements with underlying network $G$ and distinguished vertex
pair $(v,o)$. We endow this space with the topology of local
convergence in the same way as ${\mathcal{G}_*}$, except that for the
isomorphism between the local neighborhoods of two graphs, we require
that the distinguished ordered vertex pair of one graph maps to the
distinguished pair of the other graph. There is\vadjust{\goodbreak} a suitable metric for
this convergence that makes ${\mathcal{G}_{**}}$ a complete separable
metric space.

A probability measure $\mu$ on ${\mathcal{G}_*}$ is called \textit{unimodular}
if it satisfies the following for all Borel $f\dvtx {\mathcal
{G}_{**}}\rightarrow
[0,\infty]$:
\[
\int\sum_{x\in V(G)}f(G,o,x)\, \mathrm{d}\mu(G,o)=
\int \sum_{x\in
V(G)}f(G,x,o)\, \mathrm{d}\mu(G,o).
\]
A measure $\mu$ on ${\mathcal{G}_*}$ that satisfies the above for all Borel
$f$ supported on $ \{(G,x,y)|x\sim y \}$ is said to be \textit{involution
invariant}. It is clear that the set of unimodular measures is a subset
of the set of involution invariant measures. Proposition~2.2 of \cite
{AldLyo2007} shows that involution invariance is equivalent to unimodularity.

Involution invariance is characterized alternatively in \cite
{AldSte2004} as follows. Given a measure $\mu$ on ${\mathcal{G}_*}$,
define a
measure $\mu^*$ on ${\mathcal{G}_{**}}$ by letting its marginal
measure on
${\mathcal{G}_*}$ to be $\mu$ and the conditional measure on the second
vertex given a rooted geometric network $G$ to be the counting measure
on the neighbors of the root of~$G$. Specifically,
\[
\mu^*(\cdot)=\int_{\mathcal{G}_*}\sum_{v\sim o}
\mathbf{1}_{
\{(G,o,v)\in\cdot \}}\, \mathrm{d}\mu(G,o).
\]
Then $\mu$ is involution invariant if $\mu^*$ is invariant under the
\textit{involution} transformation
\[
\imath\dvtx {\mathcal{G}_{**}}\rightarrow{\mathcal{G}_{**}},
\imath (G,o,v)=(G,v,o).
\]
Involution $\imath$ swaps the order of the distinguished pair of
vertices, leaving all else unchanged.

The definitions carry forward when the graphs in ${\mathcal{G}_*}$ are
appended with maps from their edge sets to a complete separable metric
space. An edge cover $C$ on a graph $G$ can be represented as the graph
$G$ with a map on the edge set of $G\dvtx  e\mapsto\mathbf{1}_{ \{
e\in C \}}$.
We say that a random edge cover $C$ on a random graph $G$ is involution
invariant if the distribution of $G$ with the above map on its edges is
involution invariant.

\if0
Suppose $\tilde{\mu}_n$ is a probability measure on ${\mathcal{G}}$
concentrated on the elements of ${\mathcal{G}}$ with exactly $n$ vertices.
Let $G_n$ be a random element of ${\mathcal{G}}$ with distribution
$\tilde
{\mu}_n$. Pick one vertex of $G_n$ uniformly at random as root. Then
the component of the root is a random element of ${\mathcal{G}_*}$, with
distribution, say, $\mu_n$. If the sequence of such measures $\mu_n$
has a local weak limit in ${\mathcal{G}_*}$ as $n\rightarrow\infty$,
then in
\cite
{AldSte2004} the limit is said to be obtained by the \textit{standard
construction}.

Theorem~\ref{thmpwitConv} implies that the PWIT distribution is
obtained by the standard construction as the weak limit of a sequence
of randomly rooted complete graphs $\widebar{K}_n$.
\fi

In our model, the complete graphs $\widebar{K}_n$ are randomly rooted.
Write $C_n^*$ for the minimum-cost edge over on $\widebar{K}_n$ having the
same root as $\widebar{K}_n$. By symmetry it is easy to see that its
distribution is involution invariant. From Section~5.2 of \cite
{AldSte2004}, we see that involution invariance is preserved under weak
limits in the metric space ${\mathcal{G}_*}$ appended with the $
\{0,1 \}$-map on the edge set. Consequently, if the
sequence $C_n^*,n\ge1$, converges to an element $C^*$, then the
distribution of $C^*$ will be involution invariant. This motivates us
to study involution invariant edge covers on the limit PWIT.

\if0
In our model, the complete graphs $\widebar{K}_n$ are randomly rooted. Write
$C_n^*$ for the component of the root of $\widebar{K}_n$ in the
minimum-cost edge cover of $\widebar{K}_n$, with the same root as $\widebar{K}_n$.
Then $C_n^*$ is
a random element of ${\mathcal{G}_*}$, and by symmetry it is easy to
see that
its distribution is involution invariant. From Section~5.2 of \cite
{AldSte2004}, we see that involution invariance is preserved under weak
limits in the metric space ${\mathcal{G}_*}$. Consequently, if the sequence
$C_n^*,n\ge1$ converges to an element $C^*$ in ${\mathcal{G}_*}$,
then the
distribution of $C^*$ will be involution invariant. This motivates us
to study involution invariant edge covers on the limit PWIT.
\fi
%s6 #&#
%s6 ###
\section{Optimal involution invariant edge cover on the PWIT} \label{sectheEdgeCover}
%s6.1 #&#
%s6.1 ###
\subsection{A tree process based on the RDE}
In the PWIT we split each undirected edge into two directed edges. For
a general graph $G$, we use the notation $\overrightarrow{E}(G)$ to
denote the
set of directed edges so obtained. If $\xi_e$ is the cost of the
undirected edge\vadjust{\goodbreak} $e= \{v,w \}$, we assign the same cost to
both of the
corresponding directed edges and write the costs as $\xi
(u,v )=\xi (v,u )=\xi_e$. To each directed edge
$\overrightarrow{e}=(u,v)$, we
will assign a random variable denoted by $X
(\overrightarrow{e} )$ or
$X  (u,v )$. Typically, $X
(u,v )$ will be different from $X  (v,u )$.
The $X$ process is constructed in the following lemma, which is an
analogue of Lemma~5.8 of \cite{AldSte2004} and is proved similarly. We
include the proof here for completeness.

%
%le1 #&#
\begin{lemma}\label{lemtreeProcess}
There exists a process
\[
\bigl(\mathcal{T},\bigl(\xi_e,e\in E(\mathcal{T})\bigr),\bigl(X
  (\overrightarrow{e} ),\overrightarrow{e}\in\overrightarrow {E}(
\mathcal{T} )\bigr) \bigr),
\]
where $\mathcal{T}$ is a PWIT with edge lengths $ \{\xi_e,e\in
E(\mathcal{T} ) \}$, and $ \{X
(\overrightarrow{e}
),\overrightarrow{e}\in
\overrightarrow{E}(\mathcal{T}) \}$ is a stochastic process
satisfying the
following properties:
\begin{longlist}[(a)]
\item[(a)] For each directed edge $(u,v)\in\overrightarrow{E}(\mathcal{T})$,
%
%e16 #&#
%e16 ###
\begin{equation}
\label{eqXrecursion} X  (u,v )=\min \bigl\{ \bigl(\xi (v,w )-X
(v,w ) \bigr)^+:(v,w)\in \overrightarrow{E}(\mathcal{T}),w\neq u \bigr\}.
\end{equation}
\item[(b)] If $(u,v)\in\overrightarrow{E}(\mathcal{T})$ is directed away
from the root of
$\mathcal{T}$, then $X  (u,v )$ has the distribution
$F_*$ as in (\ref{eqfpdist}).
\item[(c)] If $(u,v)\in\overrightarrow{E}(\mathcal{T})$, the random
variables $X  (u,v )$
and $X  (v,u )$ are independent.
\item[(d)] For a fixed $z>0$, conditional on the event that there exists an
edge of length $z$ at the root, say $ \{\phi,v_z \}$, the random
variables $X  (\phi,v_z )$ and $X  (v_z,\phi )$ are independent random
variables, each having the distribution $F_*$.
\end{longlist}
\end{lemma}
\begin{pf}
Fix an integer $d\ge1$. We create independent random variables from
the distribution $F_*$, and assign one to each directed edge $(v,w)$ of
$\mathcal{T}$ where $v$ is at depth $d-1$, and $w$ is at depth $d$
from the
root. Then if $d>1$, use relation~(\ref{eqXrecursion}) to recursively
define random variables $X  (t,u )$, where $t\sim
u$ are vertices of
$\mathcal{T}$ within depth $d$ from the root. This generates a
collection of
random variables $\mathscr{C}_d$ whose joint distribution satisfies
properties (a), (b) and (c) in the statement of the lemma for all
vertices of $\mathcal{T}$ up to a depth $d$ from the root. It is easy
to see
that the sequence of collections $ \{\mathscr{C}_d,d\ge1 \}$
satisfies the conditions of Kolmogorov consistency theorem. So there
exists a collection $\mathscr{C}_\infty$ such that the restriction to
random variables corresponding to vertices up to depth $d$ is equal in
distribution to the collection $\mathscr{C}_d$ for each $d\ge1$. This
implies that random variables in $\mathscr{C}_\infty$ satisfy the
properties (a), (b) and (c).

To prove property (d), observe that a Poisson process conditioned to
have a point at $z$ is also a Poisson process of the same intensity
when that point is removed. Now conditional on the existence of the
edge $ \{\phi,v_z \}$ of length $z$, if we remove this edge
the PWIT
splits into two subtrees. Letting $\phi$ and $v_z$ to be the roots of
these two subtrees, we find that the two subtrees are independent
copies of the original PWIT~$\mathcal{T}$. From the construction in the
previous paragraph, it is clear that conditionally the random variables
$X  (\phi,v_z )$ and $X
(v_z,\phi )$ are independent, and have the
same distribution $F_*$.
\end{pf}

%s6.2 #&#
%s6.2 ###
\subsection{An involution invariant edge cover on the PWIT}
We use the process $ \{X  (\overrightarrow
{e} ) \}$ to construct an
edge cover $\mathcal{C}_{\mathrm{opt}}$ on $\mathcal{T}$.

For each vertex $v$ of the PWIT, define a set
%
%e17 #&#
%e17 ###
\begin{equation}
\label{eqCoptdef} \mathcal{C}_{\mathrm{opt}}(v)=\argmin_{y\sim v} \bigl\{
\bigl(\xi (v,y )-X  (v,y ) \bigr)^+ \bigr\}.
\end{equation}
In words, include in $\mathcal{C}_{\mathrm{opt}}(v)$ all $y\sim v$
such that $\xi (v,y )-X  (v,y )<0$,
and if there is no such $y$,
then $\mathcal{C}_{\mathrm{opt}}(v)= \{w \}$ where
$w$ is the unique (with probability~1) neighbor of $v$ that minimizes
$\xi (v,\cdot )-X  (v,\cdot )$.
Alternatively,
%
%e18 #&#
%e18 ###
\begin{equation}
\mathcal{C}_{\mathrm{opt}}(v)=\argmin_{A} \biggl\{\sum
_{y\in A}\bigl(\xi (v,y )-X  (v,y )\bigr)\dvtx A\subset
N_{v},A \mbox{ nonempty} \biggr\}.
\end{equation}
Define the edge cover to be
\[
\mathcal{C}_{\mathrm{opt}}=\bigcup_v \bigl\{
\{v,w \}\dvtx w\in\mathcal{C}_{\mathrm{opt}}(v) \bigr\}.
\]

The following lemma reassures us that the chosen edge cover does not
include wasteful edges.
%
%le2 #&#
\begin{lemma} \label{lemconsistentEC}
For any two vertices $v,w$ of $\mathcal{T}$, we have
\[
v\in\mathcal{C}_{\mathrm{opt}}(w)\quad\iff\quad\xi (v,w )<X  (v,w )+X
  (w,v ).
\]
As a consequence,
\[
v\in\mathcal{C}_{\mathrm{opt}}(w)\quad\iff\quad w\in\mathcal{C}_{\mathrm{opt}}(v).
\]
\end{lemma}
\begin{pf}
Suppose $w\in\mathcal{C}_{\mathrm{opt}}(v)$. If $\xi
(v,w )<X  (v,w
)$ then, since $X  (w,v )\ge0$, we have $\xi
(v,w )<X  (v,w )+X
(w,v )$.

If $\xi (v,w )\ge X  (v,w )$, then definition
(\ref{eqCoptdef}) of
$\mathcal{C}_{\mathrm{opt}}(v)$ and $w$'s membership to this set
implies that $w$ is the
only element of
\[
\argmin_{y\sim v} \bigl\{ \bigl(\xi (v,y )-X  (v,y ) \bigr)^+
\bigr\},
\]
that is,
\[
\xi (v,w )-X  (v,w )< \bigl(\xi (v,y )-X  (v,y ) \bigr)^+
\qquad\mbox{for all } {y\sim v,y\neq w}.
\]
Hence,
\begin{eqnarray*}
\xi (v,w )-X  (v,w )&<&\min \bigl\{ \bigl(\xi (v,y )-X  (v,y
) \bigr)^+\dvtx {y\sim v,y\neq w} \bigr\}
\\
&=&X(w,v),
\end{eqnarray*}
where the last equality follows from (\ref{eqXrecursion}). We have
thus established one direction of the first statement, that is,
\[
w\in\mathcal{C}_{\mathrm{opt}}(v)\quad\Longrightarrow\quad\xi (v,w )<X  (v,w
)+X  (w,v ).
\]

Conversely, suppose that $\xi (v,w )<X
(v,w
)+X  (w,v )$.
Then $X  (w,v )>\xi (v,w )-X
(v,w )$. Also $X  (w,v )\ge0$. Therefore,
\[
X  (w,v )\ge \bigl(\xi (v,w )-X  (v,w ) \bigr)^+,
\]
that is,
\[
\min_{{y\sim v,y\neq w}} \bigl(\xi (v,y )-X  (v,y ) \bigr)^+\ge
\bigl(\xi (v,w )-X  (v,w ) \bigr)^+.
\]
It follows that
\[
w\in\argmin_{y\sim v} \bigl(\xi (v,y )-X  (v,y ) \bigr)^+
\]
and hence $w\in\mathcal{C}_{\mathrm{opt}}(v)$.
Thus we have established the first statement of the lemma, which is
\[
w\in\mathcal{C}_{\mathrm{opt}}(v)\quad\iff\quad\xi (v,w )<X  (v,w )+X
  (w,v ).
\]
The condition on the right-hand side above is symmetric in $v,w$, and
hence the second statement of the lemma is proved.
\end{pf}

The following lemma asserts that the edge cover $\mathcal{C}_{\mathrm{opt}}$ satisfies
involution invariance. See Section~\ref{secinvInv} for definition.
The proof is similar to the proof of Lemma~24 of \cite{Ald2001}.
%
%le3 #&#
\begin{lemma}
$\mathcal{C}_{\mathrm{opt}}$ is involution invariant.
\end{lemma}
\begin{pf}
Given $\xi_e,X  (\overrightarrow{e}
),\overrightarrow{e}\in\overrightarrow{E}(\mathcal{T})$, the edge
cover $\mathcal{C}_{\mathrm{opt}}$ does not depend on the vertex labels
(which are strings from $\mathcal{V}$). Relation (\ref{eqXrecursion}) for
the $X$ process is also independent of the labels of the vertices. The
proof of the lemma is then complete by showing that the measure of the
$X$ process constructed in Lemma~\ref{lemtreeProcess} is involution invariant.

From the proof of Lemma~\ref{lemtreeProcess} it is clear that the
joint distribution of $X$ process is determined by the property that
for any $d>1$,
\[
\bigl\{X  (v,w )\vert v\mbox{ at depth }d-1\mbox{ from the root}, w
\mbox{ at depth }d\mbox{ from the root}\bigr\}
\]
are independent random variables with distribution $F_*$. We need to
show that this property is invariant under the involution map.

If $\phi$ is the root (first distinguished vertex) of $\mathcal{T}$, and
$u\sim\phi$ is the second distinguished vertex, then under the
involution map, $u$ becomes the root and $\phi$ the second
distinguished vertex. Write $\mathcal{T}_u$ for the subtree containing $u$
obtained by removing the edge $ \{\phi,u \}$. For an
arbitrary Borel
set $B$, define the event
\[
A:= \bigl\{\bigl(X  (v,w ),v\mbox{ at depth }d-1\mbox { from }u,w
\mbox{ at depth }d\mbox{ from }u\bigr)\in B \bigr\}.
\]

The inverse image of $A$ in the involution map is
\begin{eqnarray*}
\imath^{-1}(A)&=&\bigl\{\bigl(X  (v_1,w_1
),v_1\in\mathcal{T} _u,v_1\mbox{ at depth
}d\mbox{ from }\phi,
\\
&&\hspace*{59pt}
w_1\mbox{ at depth }d+1\mbox{ from }\phi;
\\
&&\hspace*{9pt} X  (v_2,w_2 ),v_2\in\mathcal{T}
\setminus \mathcal{T} _u,v_2\mbox{ at depth }d-2\mbox{
from }\phi,
\\
&&\hspace*{116pt}
w_2\mbox{ at depth }d-1\mbox{ from }\phi\bigr)\in B\bigr\}.
\end{eqnarray*}
Figure~\ref{figCoptInv} shows the edges involved.
It is clear that the random variables considered above are independent
with distribution $F_*$. Consequently the measure of the set $\imath
^{-1}(A)$ equals the measure of $A$. This completes the proof. Note
that we have used here the simpler notion of involution invariance
described in Section~\ref{secinvInv} rather than \textit{spatial
invariance} as used in \cite{Ald2001}.
\end{pf}

%s6.3 #&#
%s6.3 ###
\subsection{Evaluating the cost}
In the following theorem we evaluate the cost of the edge cover
$\mathcal{C}_{\mathrm{opt}}
$ on the $\mathcal{T}$. For obvious reasons, the expectation is twice the
right-hand side of (\ref{eqthmbpconv}).

%
%f4 #&#
%f4 ###
\begin{figure}%[t]

\includegraphics{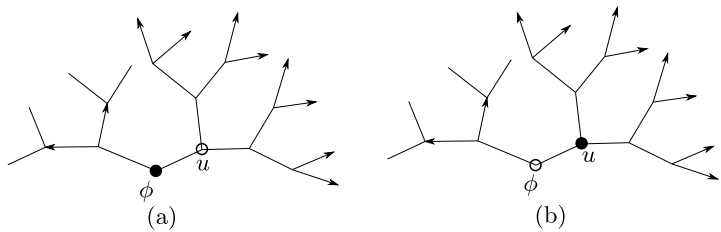}

\caption{The edges involved in events $A$ \textup{(a)} and $\imath^{-1}(A)$
\textup{(b)} are shown with arrow heads. Here $d=3$. The vertex with a filled
circle is the root, and the vertex with an unfilled circle is the
second distinguished vertex.}\label{figCoptInv}\vspace*{1pt}
\end{figure}

%
%th7 #&#
\begin{theorem}\label{thmcoptCost}
\[
\E{ \biggl[\sum_{v\in\mathcal{C}_{\mathrm{opt}}(\phi)}\xi (\phi,v )
\biggr]}=2W(1)+W(1)^2.
\]
\end{theorem}
\begin{pf}
Denote by $D$ the event that $\xi (\phi,v )>X  (\phi,v )$ for all
$v\sim\phi$. Under the event $D$, there is only one vertex in
$\mathcal{C}_{\mathrm{opt}}
(\phi)$, say $y$. By Lemma~\ref{lemconsistentEC}, $y$ is the only
neighbor of $\phi$ satisfying $\xi (\phi,y )<X  (\phi,y )+X  (y,\phi )$. Also, from (\ref
{eqXrecursion}), $X  (y,\phi )>0$. Conversely, if
there is a neighbor $y$ of $\phi$ that satisfies (i)~$X  (y,\phi )>0$,
(ii)~$\xi (\phi,y )>X  (\phi,y )$ and
(iii)~$\xi (\phi,y )<X  (\phi,y
)+X  (y,\phi )$, then from (\ref{eqXrecursion}), we have
\[
0 < X  (y,\phi ) = \min \bigl\{ \bigl(\xi (\phi,v )-X  (
\phi,v ) \bigr)^+, v \sim \phi, v \neq y \bigr\},
\]
which implies $\xi (\phi, v ) > X  (\phi,
v )$
for every $v \sim\phi, v \neq y$. This and (ii) together imply that the event $D$ holds, and
$\mathcal{C}_{\mathrm{opt}}(\phi)= \{y \}$.

Now fix a $z>0$, and condition on the event that there is a neighbor
$v_z$ of $\phi$ with $\xi (\phi,v_z )=z$. Call this event
$E_z$. If
we condition a Poisson process to have a point at some location, then
the conditional process on removing this point is again a Poisson
process with the same intensity. This shows that under $E_z$,
$X  (\phi,v_z )$ and $X
(v_z,\phi )$ both have the same distribution $F_*$. Also
they are independent. Using these facts and the characterization of the
event $D$ in the previous paragraph, the expected cost under $D$ can be
written as
%
%e19 #&#
%e19 ###
\begin{eqnarray}
\label{eqexpPart1} \qquad && \E{ \biggl[ \biggl(\sum_{v\in\mathcal{C}_{\mathrm{opt}}(\phi)}\xi_{\mathcal{T}} (\phi,v ) \biggr)\mathbf {1}_{D} \biggr]}\nonumber
\\[-1pt]
&&\qquad =\int_{z=0}^\infty zP{ \bigl\{X
(v_z,\phi )>0,z>X  (\phi,v_z ),z<X  (
\phi,v_z )+X  (v_z,\phi ) \bigr\}}\,
\mathrm{d}z\nonumber
\\[-1pt]
&&\qquad =\int_{z=0}^\infty \biggl( zP{ \bigl\{X  (\phi,v_z )=0 \bigr\}}\P{ \bigl\{X  (v_z,\phi )>z
\bigr\}}\nonumber
\\[-1pt]
&&\hspace*{24pt}\quad\qquad {} +\int
_{x=0}^z zP{ \bigl\{ X
(v_z,\phi )>z-x \bigr\}}\, \mathrm{d}F_*(x) \biggr)\,
\mathrm{d}z
\\[-1pt]
&&\qquad =\int_{z=0}^\infty \biggl(z\bigl(1-W(1)
\bigr)W(1)e^{-z}
+\int_{x=0}^z zW(1)e^{-(z-x)}W(1)e^{-x}
\, \mathrm{d}x \biggr)\, \mathrm{d}z\nonumber
\\[-1pt]
&&\qquad =W(1) \bigl(1-W(1)\bigr)+2W(1)^2\nonumber
\\[-1pt]
&&\qquad =W(1)+W(1)^2.\nonumber
\end{eqnarray}
In the second equality above, we condition on $X  (\phi,v_z )=0$ and
$X  (\phi,v_z )=x\in(0,z)$, respectively, in
the two terms of the integrand.

Under the event $D^c$, $\mathcal{C}_{\mathrm{opt}}(\phi)$ contains
all $v$ for which $\xi (\phi,v )<X  (\phi,v )$. The expected cost over
this event is given by
%
%e20 #&#
%e20 ###
\begin{eqnarray}
\label{eqexpPart2} && \E{ \biggl[ \biggl(\sum_{v\in\mathcal{C}_{\mathrm{opt}}(\phi)}
\xi (\phi,v ) \biggr)\mathbf{1}_{D^c} \biggr]}\nonumber
\\[-1pt]
&&\qquad = \E{ \biggl[\sum_{v}\xi (\phi,v )
\mathbf{1}_{ \{\xi
(\phi,v )<X  (\phi,v ) \}
} \biggr]}\nonumber
\\[-1pt]
&&\qquad =\sum_v\E{ \bigl[\xi (\phi,v )
\mathbf{1}_{ \{\xi
(\phi,v )<X  (\phi,v ) \}
} \bigr]}\nonumber
\\[-1pt]
&&\qquad=\sum_v\int_{y=0}^\infty
\P{ \bigl\{\xi (\phi,v )>y,\xi (\phi,v )<X  (\phi,v ) \bigr\} }\,
\mathrm{d}y\nonumber
\\[-1pt]
&&\qquad=\int_{y=0}^\infty\sum
_vP{ \bigl\{y<\xi (\phi,v )<X \bigr\}}\, \mathrm{d}y\nonumber
\\[-1pt]
&&\hspace*{33pt} (X\mbox{ is a $F_*$-distributed r.v. independent of the Poisson process})
\\[-1pt]
&&\qquad=\int_{y=0}^\infty \E{ \bigl[\mbox{Number of Poisson
points in }[y,X] \bigr]}\, \mathrm{d}y\nonumber
\\[-1pt]
&&\qquad=\int_{y=0}^\infty \E{ \bigl[(X-y)^+ \bigr]}
\, \mathrm{d}y\nonumber
\\[-1pt]
&&\qquad=\int_{y=0}^\infty\int_{x=y}^\infty
\widebar{F}_*(x)\, \mathrm{d}x\, \mathrm{d}y\nonumber
\\[-1pt]
&&\qquad=\int_{y=0}^\infty\int_{x=y}^\infty
W(1)e^{-x}\, \mathrm{d}x\, \mathrm{d}y\nonumber
\\[-1pt]
&&\qquad=\int_{y=0}^\infty W(1)e^{-y}\,
\mathrm{d}y\nonumber
\\[-1pt]
&&\qquad= W(1).\nonumber
\end{eqnarray}

Combining (\ref{eqexpPart1}) and (\ref{eqexpPart2}) completes the proof.
\end{pf}

In passing, we remark that $\mathcal{C}_{\mathrm{opt}}(\phi)$ is
finite almost surely.

%s6.4 #&#
%s6.4 ###
\subsection{Optimality in the class of involution invariant edge covers}
We now show that our candidate edge cover $\mathcal{C}_{\mathrm{opt}}$ has the minimum
expected cost among involution invariant edge covers on the PWIT.
%
%th8 #&#
\begin{theorem}\label{thmCopt}
Let $\mathcal{C}$ be an involution invariant edge cover of the PWIT
$\mathcal{T}
$. Write $\mathcal{C}(\phi)$ for the set of vertices of $\mathcal
{T}$ adjacent
to the root $\phi$ in $\mathcal{C}$. Then
\[
\E{ \biggl[\sum_{v\in\mathcal{C}(\phi)}\xi (\phi,v ) \biggr]}\ge \E{
\biggl[\sum_{v\in\mathcal{C}_{\mathrm{opt}}(\phi
)}\xi (\phi,v ) \biggr]}.
\]
\end{theorem}

Let us first set up some notation that will simplify the proof steps.
For each directed edge $(v,w)$ of $\mathcal{T}$, define a random variable
%Y(v,w)=\min\setbg{\sum_{y\in%A}(\len{w,y}-\X{w,y}):A\subset\nbr{w}
%
%e21 #&#
%e21 ###
\begin{equation}
\label{eqYdef} Y(v,w)=\min \biggl\{\sum_{y\in A}\bigl(
\xi (w,y )-X  (w,y )\bigr) \bigg|
\matrix{A\subset
N_{w}\setminus \{v \},
\cr
A\mbox{ nonempty}}\biggr\},
\end{equation}
where $N_{w}$ is the set of neighbors of $w$. It is easy to see that
the random variable can be written as
\[
Y(v,w)=\cases{ \displaystyle\min_{{y\sim w,y\neq v}} \bigl\{\xi (w,y )-X
(w,y ) \bigr\}c,
\vspace*{3pt}\cr
\qquad\mbox{if }\xi (w,y )-X  (w,y )
\ge0\qquad \mbox{for all } {y\sim w,y\neq v},
\vspace*{5pt}\cr
\displaystyle \sum_{{y\sim w,y\neq v}}
\bigl(\xi (w,y )-X  (w,y )\bigr)\mathbf{1}_{ \{\xi (w,y )-X
(w,y )<0 \}},
\vspace*{3pt}\cr
\qquad \mbox{otherwise.}}
\]
Note that $(Y(v,w))^+=X  (v,w )$.

Suppose that $\E{ [\sum_{v\in\mathcal{C}(\phi)}\xi (\phi,v )
]}<\infty$. Then $\mathcal{C}(\phi)$ is a finite set with
probability~1
because $ \{\xi (\phi,v ),v\sim\phi \}$ are
points of a Poisson
process of rate 1. For such an edge cover $\mathcal{C}$, define
%
%e22 #&#
%e22 ###
\begin{equation}
\label{eqAdef} A(\mathcal{C})=\sum_{v\in\mathcal{C}(\phi)}X
(\phi,v )+\max_{v\notin\mathcal{C} (\phi),v\sim\phi}Y(v,\phi).
\end{equation}
The $\max$ operation in the above equation is over an infinite number
of vertices; however, in the remark after the proof of Lemma~\ref
{lemD}, we will show that effectively $Y(v,\phi)$ assumes only
finitely many values as we vary $v$, and hence the $\max$ operation as
well as $A(\mathcal{C})$ are almost surely well defined.

The following two lemmas will be used to prove Theorem~\ref{thmCopt}.

%
%le4 #&#
\begin{lemma}\label{lemD}
Let $\mathcal{C}$ be an edge cover rule on the PWIT such that
\[
\E{ \biggl[\sum_{v\in\mathcal{C}(\phi)}\xi (\phi,v ) \biggr]}<\infty.
\]
Then almost surely,
\[
\sum_{v\in\mathcal{C}(\phi)}\xi (\phi,v )\ge A(\mathcal{C}).
\]
Furthermore,
\[
\sum_{v\in\mathcal{C}_{\mathrm{opt}}(\phi)}\xi (\phi,v )=A(\mathcal{C}_{\mathrm{opt}}).
\]
\end{lemma}

%
%le5 #&#
\begin{lemma}\label{lemA}
Let $\mathcal{C}$ be an edge cover rule on the PWIT such that
\[
\E{ \biggl[\sum_{v\in\mathcal{C}(\phi)}\xi (\phi,v ) \biggr]} <
\infty.
\]
If $\mathcal{C}$ is involution invariant, we have
$\E{ [A(\mathcal{C}) ]}\ge \E{ [A(\mathcal{C}_{\mathrm{opt}}) ]}$.
\end{lemma}

\begin{pf*}{Proof of Theorem~\ref{thmCopt}}
If $\E{ [\sum_{v\in\mathcal{C}(\phi)}\xi (\phi,v
)
]}=\infty$,
the statement of the theorem is trivially true. Assume that it is
finite. We are now in a position to apply Lemmas \ref{lemD} and \ref
{lemA} as follows to get the result
\begin{eqnarray*}
\E{ \biggl[\sum_{v\in\mathcal{C}(\phi)}\xi (\phi,v ) \biggr]} &\ge& \E{
\bigl[A(\mathcal{C}) \bigr]}\qquad\mbox{(Lemma~\ref{lemD})}
\\
&\ge& \E{ \bigl[A(\mathcal{C}_{\mathrm{opt}}) \bigr]}\qquad\mbox
{(Lemma~\ref{lemA})}
\\
&=&\E{ \biggl[\sum_{v\in\mathcal{C}_{\mathrm{opt}}(\phi)}\xi (\phi,v ) \biggr]}\qquad
\mbox {(Lemma~\ref{lemD})}.
\end{eqnarray*}\upqed
\end{pf*}

Let us now complete the proofs of Lemmas \ref{lemD} and \ref{lemA}.
\begin{pf*}{Proof of Lemma~\ref{lemD}}
From (\ref{eqYdef}), we have
\[
Y(v,\phi)\le\sum_{y\in A}\bigl(\xi (\phi,y )-X
  (\phi,y )\bigr)
\]
for all $A\subset N_{\phi}\setminus \{v \}$, $A$ nonempty.

For any $v\notin\mathcal{C}(\phi)$, we can choose $A=\mathcal
{C}(\phi)$ to obtain
\[
Y(v,\phi)\le\sum_{y\in\mathcal{C}(\phi)}\bigl(\xi (\phi,y )-X
  (\phi,y )\bigr).
\]
This implies
%
%e23 #&#
%e23 ###
\begin{equation}
\label{eqmaxY} \max_{v\notin\mathcal{C}(\phi),v\sim\phi} Y(v,\phi)\le\sum
_{y\in\mathcal{C}(\phi)}\bigl(\xi (\phi,y )-X  (\phi,y )\bigr).
\end{equation}
Thanks to the finite expectation assumption in the lemma, $\mathcal
{C}(\phi
)$ is a finite set almost surely, and so $\sum_{y \in\mathcal
{C}(\phi)}
X  (\phi,y )$ is finite. Rearrangement of (\ref
{eqmaxY}) then yields
\[
\sum_{v\in\mathcal{C}(\phi)}\xi (\phi,v )\ge A(\mathcal{C}).
\]

Now recall the alternate characterization of $\mathcal{C}_{\mathrm{opt}}$ via
%
%e24 #&#
%e24 ###
\begin{equation}
\label{eqCoptalt} \qquad\mathcal{C}_{\mathrm{opt}}(w)=\argmin_A \biggl\{
\sum_{y\in A}\bigl(\xi (w,y )-X  (w,y )
\bigr)\dvtx A\subset N_{w},A\mbox{ nonempty} \biggr\}.
\end{equation}
From (\ref{eqYdef}) and (\ref{eqCoptalt}), for any $v\notin
\mathcal{C}_{\mathrm{opt}}
(\phi)$, we have
%
%e25 #&#
%e25 ###
\begin{equation}
\label{eqYvphivalues} Y(v,\phi)=\sum_{y\in\mathcal{C}_{\mathrm{opt}}(\phi)}\bigl(\xi (\phi,y
)-X  (\phi,y )\bigr)
\end{equation}
and hence
\[
\max_{v\notin\mathcal{C}_{\mathrm{opt}}(\phi),v\sim\phi}Y(v,\phi )=\sum_{y\in\mathcal{C}_{\mathrm{opt}}(\phi)}
\bigl(\xi (\phi,y )-X  (\phi,y )\bigr).
\]
It follows by rearrangement that
\[
\sum_{v\in\mathcal{C}_{\mathrm{opt}}(\phi)}\xi (\phi,v )=A(\mathcal{C}_{\mathrm{opt}}).
\]
\upqed
\end{pf*}

Let us quickly reassure the reader that the max operation in (\ref
{eqAdef}) is well defined. Notice that (\ref{eqYvphivalues}) implies
that $Y(w,\phi)$ takes values in the finite set
\[
\biggl\{\sum_{y\in\mathcal{C}_{\mathrm{opt}}(\phi)}\bigl(\xi (\phi,y )-X
(\phi,y )\bigr) \biggr\}\cup \bigl\{Y(v,\phi)\vert v\in\mathcal{C}_{\mathrm{opt}}(
\phi) \bigr\}.
\]
That $\mathcal{C}_{\mathrm{opt}}(\phi)$ is finite (almost surely)
can be gleaned from
Theorem~\ref{thmcoptCost}.
This validates the assertion that the $\max$ in the definition of
$A(\mathcal{C})$ is well defined.

\begin{pf*}{Proof of Lemma~\ref{lemA}}
Define
%
%e26 #&#
%e26 ###
\begin{equation}
\label{eqAtilde} \widetilde{A}(\mathcal{C})=\sum_{v\in\mathcal{C}(\phi)}X
  (v,\phi )+\max_{v\notin\mathcal{C} (\phi),v\sim\phi
}Y(v,\phi).
\end{equation}
We will prove Lemma~\ref{lemA} by showing the following two results:
\begin{longlist}[(a)]
\item[(a)] For an involution invariant edge cover $\mathcal{C}$,
%
%e27 #&#
%e27 ###
\begin{equation}
\label{eqlemA1} \E{ \bigl[\widetilde{A}(\mathcal{C}) \bigr]}=\E{ \bigl[A(\mathcal{C})
\bigr]}.
\end{equation}
\item[(b)] Almost surely,
%
%e28 #&#
%e28 ###
\begin{equation}
\label{eqlemA2} \widetilde{A}(\mathcal{C})\ge\widetilde{A}(\mathcal{C}_{\mathrm{opt}}).
\end{equation}
\end{longlist}

We first prove (\ref{eqlemA1}). First, by involution invariance of
$\mathcal{C}$, we have
%
%e29 #&#
%e29 ###
\begin{equation}
\label{eqXsymm} \E{ \biggl[\sum_{v\in\mathcal{C}(\phi)}X  (
\phi,v ) \biggr]}=\E{ \biggl[\sum_{v\in\mathcal{C}(\phi
)}X
(v,\phi ) \biggr]}.
\end{equation}
Indeed, the left-hand side equals
\[
\int_{\mathcal{G}_*}\sum_{v\sim\phi}X
  (\phi,v )\mathbf{1}_{ \{ \{\phi,v \}\in\mathcal
{C} \}}\, \mathrm{d}
\mu_\mathcal{C}\bigl([G,\phi]\bigr),
\]
where $\mu_\mathcal{C}$ is the probability measure on ${\mathcal{G}_*}$
corresponding to $(\mathcal{T},\mathcal{C})$.
By involution invariance, this equals
\[
\int_{\mathcal{G}_*}\sum_{v\sim\phi}X
  (v,\phi )\mathbf{1}_{ \{ \{v,\phi \}\in\mathcal
{C} \}}\, \mathrm{d}
\mu_\mathcal{C}\bigl([G,\phi]\bigr),
\]
which is equal to the right-hand side of (\ref{eqXsymm}).
Thanks to the finite expectation assumption of the lemma, we saw in the
proof of Lemma~\ref{lemD} that
\[
\max_{v\notin\mathcal{C}(\phi
),v\sim\phi}Y(v,\phi)
\]
is
finite almost surely. Now observe that $A(\mathcal{C})$ [resp., $\widetilde{A}(\mathcal{C})$] is obtained by adding the almost surely finite random
variable $\max_{v\notin\mathcal{C}(\phi),v\sim\phi}Y(v,\phi)$ to
the random variable which
is the argument of the expectation on the left-hand side of (\ref
{eqXsymm}) [resp., the right-hand side of (\ref{eqXsymm})]. Taking
expectation and using the equality in (\ref{eqXsymm}), we get (\ref
{eqlemA1}).

Now we will prove (\ref{eqlemA2}). First condition on the event
$L_1= \{\llvert \mathcal{C}_{\mathrm{opt}}(\phi)\rrvert >1
\}$. Observe that, under $L_1$, $\xi (\phi,y )-X  (\phi,y )<0$, $y\sim\phi$ if and
only if $y\in\mathcal{C}_{\mathrm{opt}}
(\phi)$, and there are at least two such $y$. Then, by (\ref{eqXrecursion}),
%
%e30 #&#
%e30 ###
\begin{equation}
\label{eqApart1} X  (v,\phi )=0\qquad\mbox{for all }v\sim\phi.
\end{equation}
Also, from (\ref{eqYdef}) and (\ref{eqCoptalt}),
\[
Y(v,\phi)\ge\sum_{y\in\mathcal{C}_{\mathrm{opt}}(\phi)}\bigl(\xi (\phi,y )-X
  (\phi,y )\bigr)=Y(w,\phi)
\]
if $w\notin\mathcal{C}_{\mathrm{opt}}(\phi)$. This implies
\[
Y(v,\phi)\ge\max_{w\notin\mathcal{C}_{\mathrm{opt}}(\phi),w\sim
\phi}Y(w,\phi)\qquad\mbox{for all }v\sim\phi.
\]
In particular,
%
%e31 #&#
%e31 ###
\begin{equation}
\label{eqApart2} \max_{v\notin\mathcal{C}(\phi),v\sim\phi}Y(v,\phi)\ge\max
_{w\notin\mathcal{C}_{\mathrm{opt}}(\phi),w\sim\phi}Y(w,\phi).
\end{equation}
Combining (\ref{eqApart1}) and (\ref{eqApart2}) gives
%
%e32 #&#
\begin{eqnarray}
\label{eqApart3} && \sum_{v\in\mathcal{C}(\phi)}X  (v,\phi )+
\max_{v\notin\mathcal{C}(\phi),v\sim\phi}Y(v,\phi)
\nonumber\\[-8pt]\\[-8pt]
&&\qquad \ge \sum_{v\in\mathcal{C}_{\mathrm{opt}}(\phi)}X  (v,\phi )+\max
_{v\notin\mathcal{C}_{\mathrm{opt}}(\phi
),v\sim\phi}Y(v,\phi).\nonumber
\end{eqnarray}
Thus $\widetilde{A}(\mathcal{C})\ge\widetilde{A}(\mathcal{C}_{\mathrm{opt}})$ under $L_1$.

Now consider the event $L_2= \{\llvert \mathcal{C}_{\mathrm{opt}}(\phi)\rrvert =1 \}$. Let
\begin{eqnarray*}
X_\phi^{(1)}&=&\min_{v\sim\phi}\bigl(\xi (\phi,v
)-X  (\phi,v )\bigr)\quad\mbox{and}
\\
X_\phi^{(2)}&=&{\min_{v\sim\phi}}^{(2)}
\bigl(\xi (\phi,v )-X  (\phi,v )\bigr),
\end{eqnarray*}
where ${\min_{}}^{(2)}$ stands for the second minimum.

Let $\mathcal{C}_{\mathrm{opt}}(\phi)= \{u \}$. Then
$X  (u,\phi
)=X_\phi^{(2)}$, and for $v\in
\mathcal{C}(\phi)\setminus\mathcal{C}_{\mathrm{opt}}(\phi)$,
$X  (v,\phi
)=(X_\phi^{(1)})^+$. So
we get
%
%e33 #&#
\begin{eqnarray}
\label{eqApart4} && \sum_{v\in\mathcal{C}(\phi)}X  (v,\phi )-\sum
_{v\in\mathcal{C}_{\mathrm{opt}}(\phi)}X  (v,\phi )
\nonumber\\[-8pt]\\[-8pt]
&&\qquad =\sum_{v\in\mathcal{C}(\phi)\setminus\mathcal{C}_{\mathrm{opt}}(\phi)}\bigl(X_\phi ^{(1)}
\bigr)^+-X_\phi^{(2)} \mathbf{1}_{ \{u\notin\mathcal{C}(\phi) \}}.\nonumber
\end{eqnarray}

If $v\notin\mathcal{C}_{\mathrm{opt}}(\phi)$, then $Y(v,\phi
)=X_\phi^{(1)}$. Also
$Y(u,\phi
)=X_\phi^{(2)}$. Since $X_\phi^{(2)}\ge X_\phi^{(1)}$, we get
\[
\max_{v\notin\mathcal{C}(\phi),v\sim\phi}Y(v,\phi)=X_\phi ^{(2)}
\mathbf{1}_{ \{u\notin\mathcal{C} (\phi) \}
}+X_\phi^{(1)}\mathbf{1}_{ \{u\in\mathcal{C}(\phi) \}}
\]
and
\[
\max_{v\notin\mathcal{C}_{\mathrm{opt}}(\phi),v\sim\phi}Y(v,\phi )=X_\phi^{(1)}.
\]
Therefore,
%
%e34 #&#
%e32 ###
\begin{equation}
\label{eqApart5} \max_{v\notin\mathcal{C}(\phi),v\sim\phi}Y(v,\phi)-\max_{v\notin\mathcal{C}_{\mathrm{opt}}(\phi),v\sim\phi}Y(v,
\phi )= \bigl(X_\phi^{(2)}-X_\phi^{(1)}
\bigr)\mathbf{1}_{ \{u\notin
\mathcal{C}(\phi) \}}.
\end{equation}
Adding (\ref{eqApart4}) and (\ref{eqApart5}), and canceling $X_\phi^{(2)}
\mathbf{1}_{ \{u\notin\mathcal{C}(\phi) \}}$, we get
\begin{eqnarray*}
&& \sum_{v\in\mathcal{C}(\phi)}X  (v,\phi )+\max
_{v\notin\mathcal{C}(\phi),v\sim\phi}Y(v,\phi)
-\sum_{v\in\mathcal{C}_{\mathrm{opt}}(\phi)}X  (v,\phi )-\max
_{v\notin\mathcal{C}_{\mathrm{opt}}(\phi
),v\sim\phi}Y(v,\phi)
\\
&&\qquad =\sum_{v\in\mathcal{C}(\phi)\setminus\mathcal{C}_{\mathrm{opt}}(\phi)}\bigl(X_\phi ^{(1)}
\bigr)^+-X_\phi^{(1)} \mathbf{1}_{ \{u\notin\mathcal{C}(\phi) \}}
\ge 0,
\end{eqnarray*}
where the last inequality follows because there exists a
$v \in\mathcal{C}(\phi) \setminus\mathcal{C}_{\mathrm{opt}}(\phi)$
by virtue of our assumption that
$\mathcal{C}(\phi) \neq\mathcal{C}_{\mathrm{opt}}(\phi)$.
Thus $\widetilde{A}(\mathcal{C})\ge\widetilde{A}(\mathcal{C}_{\mathrm{opt}})$ under $L_2$ as well.
\end{pf*}

%s7 #&#
%s7 ###
\section{Completing the lower bound}\label{seclowerBd}
In the previous section we described an edge cover $\mathcal
{C}_{\mathrm{opt}}$ on the
infinite tree $\mathcal{T}$. We showed that this edge cover satisfies the
expected property of involution invariance, and it has the minimum
expected cost among all edge covers having this property. We use this
to show now that the expected cost of $\mathcal{C}_{\mathrm{opt}}$
serves as an asymptotic
lower bound on the expected cost of min-cost edge covers on $\widebar{K}_n$.

%
%th9 #&#
\begin{theorem}\label{thmlowerBd}
Let $C_n^*$ be the optimal edge cover on $\widebar{K}_n$. Then
\[
\liminf_{n\rightarrow\infty}\E{ \biggl[\sum_{ \{\phi,v \}
\in C_n^*}
\xi_{\widebar{K}_n} (\phi,v ) \biggr]}\ge2W(1)+W(1)^2.
\]
\end{theorem}
\begin{pf}
Take a subsequence $ \{n_k,k\ge1 \}$ for which the $\liminf
$ above
is a limit. Now consider the joint sequence $(C_{n_k}^*,\widebar{K}
_{n_k})_{k\ge1}$ in ${\mathcal{G}_*}\times{\mathcal{G}_*}$. Because
$\widebar{K}
_{n_k}\xrightarrow{\mathrm{l.w.}}\mathcal{T}$, for every $\varepsilon
>0$ there is a compact subset
$\mathcal{K}$ of ${\mathcal{G}_*}$, with $\P \{\widebar
{K}_{n_k}\in\mathcal{K} \}>1-\varepsilon$ for all $k$. Also, we
can take the graphs $\widebar{K}
_{n_k}$ to be on a common vertex set $\widetilde{\mathcal{V}}$, and
assume that all graphs in $\mathcal{K}$ are defined on the same vertex
set. Let $\widetilde{\mathcal{E}}$ denote the set of all possible edges.
Let $\mathcal{K}_S$ denote the set $ \{H\mbox{ is a subgraph of
}G\vert G\in\mathcal{K} \}$. Since $C_{n_k}^*$ is a subgraph of
$\widebar{K}
_{n_k}$, $\P \{C_{n_k}^*\in\mathcal{K}_S \}>1-\varepsilon$
for all
$k$. An element of $\mathcal{K}_S$ can be identified with an element
of $\mathcal{K}\times \{0,1 \}^{\widetilde{\mathcal{E}}}$,
where 1 or
0 denotes the presence or absence of an edge, respectively. Since the
latter is a compact set, so is $\mathcal{K}_S$. This shows that the
sequence of random graphs $ \{C_{n_k}^* \}_{k\ge1}$ is
tight. By
completeness of ${\mathcal{G}_*}$, we have that $ \{
(C_{n_k}^*,\widebar{K} _{n_k}),k\ge1 \}$ is sequentially
compact. Therefore, there exists a
further subsequence $ \{n_j,j\ge1 \}$ of $ \{n_k,k\ge
1 \}$ such
that $(C_{n_j}^*,\widebar{K}_{n_j})$ converges in the local weak sense to
$(C^*,\mathcal{T})$. Since the $C_n^*$ distribution is involution invariant,
so is the distribution of $C^*$. By Skorohod's theorem we can assume
the convergence occurs almost surely in some probability space. By the
definition of local weak convergence
\[
\sum_{ \{\phi,v \}\in C_{n_j}^*}\xi_{\widebar{K}
_{n_j}} (\phi,v )\rightarrow
\sum_{v\in C^*(\phi)}\xi_{\mathcal{T}} (\phi,v )\qquad\mbox
{as }n\rightarrow\infty \mbox{ a.s.}
\]
By Fatou's lemma
\[
\liminf_{j\rightarrow\infty}\E{ \biggl[\sum_{ \{\phi,v \}
\in C_{n_j}^*}
\xi_{\widebar{K}_{n_j}} (\phi,v ) \biggr]}\ge \E{ \biggl[\sum
_{v\in C^*(\phi)}\xi_{\mathcal{T}} (\phi,v ) \biggr]}.
\]
By Theorems~\ref{thmCopt}~and~\ref{thmcoptCost},
\[
\E{ \biggl[\sum_{v\in C^*(\phi)}\xi_{\mathcal{T}} (\phi,v )
\biggr]}\ge \E{ \biggl[\sum_{v\in\mathcal{C}_{\mathrm{opt}}(\phi)}\xi_{\mathcal{T}} (
\phi,v ) \biggr]} =2W(1)+W(1)^2.
\]
This completes the proof.
\end{pf}

%%\input{theMatching}
%%\input{the2factor}
%s8 #&#
%s8 ###
\section{Belief propagation} \label{secbp}
To prove the upper bound on $EC_n$ in order to complete the proof
of Theorem~\ref{thmlimit}, we will construct edge covers on $K_n,n\ge1$, with costs $W(1)+W(1)^2/2+o(1)$. This is achieved using
belief propagation as described in Section~\ref{secresults}.

We follow the approach of \cite{SalSha2009} to prove Theorem~\ref
{thmbpconv}. In this section we will show the convergence of the BP
algorithm on the PWIT $\mathcal{T}$, and relate the converged solution with
the edge cover $\mathcal{C}_{\mathrm{opt}}$ of Section~\ref
{sectheEdgeCover}. In the next
section we show that the belief propagation on $\widebar{K}_n$
converges to
belief propagation on $\mathcal{T}$ as $n\rightarrow\infty$.
%s8.1 #&#
%s8.1 ###
\subsection{Convergence of BP on the PWIT}
In this section we will prove that the messages on $\mathcal{T}$ converge,
and relate the resulting edge cover with the cover $\mathcal
{C}_{\mathrm{opt}}$ of
Section~\ref{sectheEdgeCover}.

The message process can essentially be written as
%
%e35 #&#
%e33 ###
\begin{equation}
\label{eqBPonT} X_{\mathcal{T}}^{k+1} (\dot{v},v )=\min
_{i\ge1} \bigl\{ \bigl(\xi_{\mathcal{T}} (v,v.i
)-X_{\mathcal
{T}}^{k} (v,v.i ) \bigr)^+ \bigr\},
\end{equation}
where the initial messages $X_{\mathcal{T}}^{0} (\dot{v},v )$
are i.i.d. random
variables [zero in the case of our algorithm; see (\ref{eqbpInit})].

By the structure of $\mathcal{T}$, it is clear that for a fixed $k\ge0$,
all the messages $X_{\mathcal{T}}^{k} (\dot{v},v ),v\in
\mathcal
{V}$ share the same
distribution. Also, it can be seen from the analysis of RDE (\ref
{eqECrde}) in Section~\ref{secrde} that if we denote the
complementary c.d.f. of this distribution at some step $k$ by $\widebar{F}$,
then after one update the complementary c.d.f. is given by the map
\[
T\widebar{F}(y)=\cases{ \displaystyle e^{-y}\exp \biggl(-\int
_0^\infty\widebar{F}(t)\, \mathrm{d}t \biggr), &
\quad if $y\ge0$,
\cr
1, &\quad if $y<0$.}
\]
The operator $T$ thus defined on the space $\mathcal{D}$ of
complementary c.d.f.'s of $\widebar{\mathbf{R}}$-valued random variables has a
unique fixed point $\widebar{F}_*$ given by (\ref{eqfpdist}).

The following theorem shows that the fixed point $\widebar{F}_*$ has the
full space $\mathcal{D}$ as its domain of attraction. In other words,
irrespective of the initial distribution, the common distribution of
the messages $X_{\mathcal{T}}^{k} (\dot{v},v ),v\in
\mathcal
{V}$ converges to the
distribution $F_*$ as $k\rightarrow\infty$.

%
%th10 #&#
\begin{theorem} \label{thmfpConv}
For any $\widebar{F}\in\mathcal{D}$,
\[
\lim_{k\rightarrow\infty}T^k\widebar{F}=\widebar{F}_*.
\]
\end{theorem}
\begin{pf}
For any $y\ge0$ and $k\ge0$,
\[
T^{k+1}\widebar{F}(y)=e^{-y}\exp \biggl(-\int
_0^\infty T^{k}\widebar{F}(t)\,
\mathrm{d}t \biggr).
\]
Thus for $k\ge1$, $T^k\widebar{F}(y)=c_k e^{-y}$, where $c_k,k\ge1$, are
nonnegative real numbers satisfying
\[
c_{k+1}=\exp \biggl(-\int_0^\infty
c_k e^{-t}\, \mathrm{d}t \biggr)=e^{-c_k}.
\]
It is easy to check that $c_k\rightarrow W(1)$. Consequently, $T^k\widebar{F}\rightarrow
\widebar{F}_*$.
\end{pf}
%
%s8.2 #&#
%s8.2 ###
\subsection{Endogeny and bivariate uniqueness}
We have established the convergence of the messages on $\mathcal{T}$ in
distribution. We now ask for the joint convergence of the message
process on the tree. In particular, the question is whether there is a
limit process satisfying the requirements of Lemma~\ref{lemtreeProcess}.

An important property of the limiting process that allows us to come to
this conclusion is \textit{endogeny} introduced in \cite{AldBan2005}.
Endogeny is a property of the recursive tree process (RTP) that it is
measurable with respect to the i.i.d. process $(\xi_{\bi i},N_{\bi
i}),i\in\mathcal{V}$.

\begin{defnn*}
An invariant RTP with marginal distribution $\mu$ is said to be
endogenous if the root variable $X_\phi$ is almost surely measurable
with respect to the $\sigma$-algebra
\[
\sigma \bigl( \bigl\{(\xi_{\bi i},N_{\bi i})|{\bi i}\in\mathcal
{V} \bigr\} \bigr).
\]
\end{defnn*}
Endogeny is related to another property of the RTP termed as \textit{bivariate uniqueness} again introduced in \cite{AldBan2005}.

For a general RDE (\ref{eqrde}) write $T\dvtx \mathcal{P}\rightarrow
\mathcal
{P}(S)$ for the map induced by the function~$g$. Let $\mathcal
{P}^{(2)}$ denote the\vspace*{2pt} space of probability measures on $S\times S$ with
marginals in $\mathcal{P}$. We now define a bivariate map
$T^{(2)}\dvtx \mathcal{P}^{(2)}\rightarrow\mathcal{P}(S\times S)$, which
maps a
distribution $\mu^{(2)}\in\mathcal{P}^{(2)}$ to the joint
distribution of
\[
\pmatrix{g\bigl(\xi;\bigl(X^{(1)}_j,1\le
j< N\bigr)\bigr)
\vspace*{3pt}\cr
g\bigl(\xi;\bigl(X^{(2)}_j,1\le j< N\bigr)\bigr)},
\]
where\vspace*{-1pt} $(X^{(1)}_j,X^{(2)}_j)_{j\ge1}$ are independent with joint
distribution $\mu^{(2)}$ on $S\times S$, and the family of random
variables $(X^{(1)}_j,X^{(2)}_j)_{j\ge1}$ are independent\vspace*{1pt} of the pair
$(\xi,N)$.\vadjust{\goodbreak}

It is easy to see that if $\mu$ is a fixed point of the RDE, then the\vspace*{1pt}
associated \textit{diagonal measure} $\mu^\nearrow:=\operatorname{Law}(X,X)$
where $X\sim\mu$ is a fixed point of the operator $T^{(2)}$.

\begin{defnn*}
An invariant RTP with marginal distribution $\mu$ is said to have the
bivariate uniqueness property if $\mu^\nearrow$ is the unique fixed point
of the operator $T^{(2)}$ with marginals $\mu$.
\end{defnn*}

Theorem~11 of \cite{AldBan2005} stated below shows that under certain
assumptions, endogeny and bivariate uniqueness are equivalent.

%
%th11 #&#
\begin{theorem}[(Theorem 11 of \cite{AldBan2005})]\label{thmendBVU}
Let $S$ be a Polish space. Consider an invariant RTP with
marginal distribution $\mu$:

\begin{longlist}[(a)]
\item[(a)] If the endogenous property holds, then the bivariate uniqueness
property holds.
\item[(b)] Conversely, suppose the bivariate uniqueness property holds. If
also $T^{(2)}$ is continuous with respect to weak convergence on the
set of bivariate distributions with marginals $\mu$, then the
endogenous property holds.
\item[(c)] The endogenous property holds if and only if ${T^{(2)}}^n(\mu
\otimes\mu)\xrightarrow{\mathrm{D}}\mu^\nearrow$, where $\mu
\otimes\mu$ is the
product measure.
\end{longlist}
\end{theorem}

The following theorem establishes the endogeny of the edge cover RDE.
%
%th12 #&#
\begin{theorem}\label{thmendOfECrde}
The invariant RTP with marginal $\mu_*$ (with c.d.f. $F_*$) associated
with the edge cover RDE (\ref{eqECrde}) is endogenous.
\end{theorem}
\begin{pf}
By Theorem~\ref{thmendBVU}(b) it is sufficient to prove bivariate
uniqueness and continuity for the map $T^{(2)}\dvtx \mathcal{P}(\mathbf{R}_{+}
\times\mathbf{R}_{+})\rightarrow\mathcal{P}(\mathbf{R}_{+}\times
\mathbf{R}_{+})$, where
$\mathbf{R}_{+}
=[0,\infty)$ and $T^{(2)}(\mu^{(2)})$ is the distribution of
\[
\pmatrix{ X
\cr
Y}
= %
\pmatrix{ \displaystyle\min
_{i\ge1}(\xi_i-X_i)^+
\vspace*{3pt}\cr
\displaystyle\min_{i\ge1}(\xi_i-Y_i)^+ },
\]
where $(X_i,Y_i)_{i\ge1}$ are independent with joint distribution $\mu
^{(2)}$ on $\mathbf{R}_{+}^2$, and are independent of $(\xi_i)_{i\ge1}$
which are points of a Poisson process of rate 1 on~$\mathbf{R}_{+}$.

To prove bivariate uniqueness, we have to show that if $\mu_*^{(2)}$
is a fixed point of the above map (with marginals $\mu_*$), then $X=Y
\mbox{ a.s. } (\mu_*^{(2)})$. By Lemma~1 of~\cite{Ban2011} this is
equivalent to showing $X\eqdist Y\eqdist X\wedge Y$. Let
$(X_i,Y_i)_{i\ge1}$ be i.i.d. with distribution $\mu^{(2)}$. The set
of points $\mathcal{P}:= \{(\xi_i;(X_i,Y_i))|i\ge1 \}$
forms a Poisson\vspace*{2pt} process on $(0,\infty)\times\mathbf{R}_{+}^2$ with intensity
$\, \mathrm{d}t \mu_*^{(2)}(\, \mathrm{d}(x,y))$
at $(t;(x,y))$. Writing $G(x,y)=\P{ \{X>x,Y>y \}}$ for $x,y\in\mathbf{R}_{+}$, we get
%
%e36 #&#
\begin{eqnarray}
\label{eqGxy} G(x,y)&=&\P{ \{\xi_i-X_i>x,
\xi_i-Y_i>y,\mbox{ for all }i\ge 1 \}}\nonumber
\\
&=&\P{ \bigl\{\mbox{No point of $\mathcal{P}$ in } \bigl\{ \bigl(t;(u,v)\bigr)\dvtx t-u
\le x\mbox{ or }t-v\le y \bigr\} \bigr\}}\nonumber
\\
&=&\exp \biggl(-\int_{t=0}^{x\vee y}
\mathrm{d}t-\int_{t=x\vee y}^\infty \P{ \{t-X_1
\le x\mbox{ or }t-Y_1\le y \}}\, \mathrm{d}t \biggr)\nonumber
\\
&=&e^{-x\vee y}\exp \biggl(-\int_{t=x\vee y}^\infty \P{
\{X_1\ge t-x\mbox{ or }Y_1\ge t-y \}}\,
\mathrm{d}t \biggr)
\nonumber\\[-8pt]\\[-8pt]
&=&e^{-x\vee y}\exp \biggl(-\int_{t=x\vee y}^\infty
\bigl(W(1)e^{-(t-x)}+W(1)e^{-(t-y)}\nonumber
\\
&&\hspace*{94pt} {}-\P{ \{X_1\ge t-x,Y_1\ge t-y \}} \bigr)\,
\mathrm{d}t \biggr)\nonumber
\\
&=&e^{-x\vee y}\exp\bigl(-W(1)e^{-x\vee y}\bigl(e^x+e^y
\bigr)\bigr)\nonumber
\\
&&{}\times \exp \biggl(\int_{t=x\vee y}^\infty \P{
\{X_1\ge t-x,Y_1\ge t-y \}}\, \mathrm{d}t
\biggr).\nonumber
\end{eqnarray}
From this, setting $x=y$, it is clear that $G(x,x) = c e^{-x}, x \ge
0$, for some constant~$c$. We now have to evaluate the constant.

Observe that the only place where $G(x,x)$ can be discontinuous (if at
all) is at $x = 0$. As a consequence, with $x = y$ and the change of
variable $z = t - x$, we see that the integral inside the exponent in
(\ref{eqGxy}) is $\int_0^{\infty} P(X_1 \ge z, Y_1 \ge z) \, \mathrm{d}z =
\int_0^{\infty} P(X_1 > z, Y_1 > z) \, \mathrm{d}z = \int_0^{\infty} G(z,z)
\, \mathrm{d}z$. With $x=y$ in (\ref{eqGxy}), and
integrating, we find that
\[
c=e^{-2W(1)}e^c,
\]
that is,
\[
ce^{-c}=e^{-2W(1)}.
\]
Since $W(1)=e^{-W(1)}$, it can be seen that $c=W(1)$ solves the above
equation. Because $G(0,0)\le1$, we have $c\le1$, and noting that the
function $x\mapsto xe^{-x}$ is monotone increasing for $0\le x\le1$,
we conclude that $c=W(1)$ is the only solution. Thus $G=\widebar{F}_*$,
that is, $X\wedge Y\eqdist X\eqdist Y$. This establishes bivariate uniqueness.

Now to establish endogeny it remains to prove the continuity hypothesis
of Theorem~\ref{thmendBVU}(b). Note that we require continuity of the
map $T^{(2)}$ only over the subset $\mathcal{P}_*\subset\mathcal
{P}(\mathbf{R}_{+}^2)$ which contains probability distributions with both
marginals equal to~$\mu_*$. We need to show that for any $\mu
^{(2)}\in\mathcal{P}_*$ and a sequence $(\mu_n^{(2)})_{n\ge1}$ in
$\mathcal{P}_*$ such that $\mu_n^{(2)}\xrightarrow{\mathrm{D}}\mu
^{(2)}$, we have
$T^{(2)}(\mu_n^{(2)})\xrightarrow{\mathrm{D}}T^{(2)}(\mu^{(2)})$.

Take a probability space $(\Omega,\mathcal{F},P)$ in which there are
random vectors $(X,Y)\sim\mu^{(2)}$ and a sequence of random vectors
$ \{(X_n,Y_n),n\ge1 \}$, with $(X_n,Y_n)\sim\mu_n^{(2)}$. Then
$(X_n,Y_n)\xrightarrow{\mathrm{D}}(X,Y)$. By following the steps of
(\ref{eqGxy}),
for $x,y\in\mathbf{R}_{+}$, we can write
%
%e37 #&#
%e34 ###
\begin{eqnarray}
G_n(x,y)&=& T^{(2)}\bigl(\mu_n^{2}
\bigr) \bigl((x,\infty),(y,\infty)\bigr)\nonumber
\\
&=&e^{-x\vee y}\exp\bigl(-W(1)e^{-x\vee y}\bigl(e^x+e^y
\bigr)\bigr)\nonumber
\\
&&{}\times \exp \biggl(\int_{t=x\vee y}^\infty \P{
\{X_n\ge t-x,Y_n\ge t-y \}}\, \mathrm{d}t \biggr)\nonumber
\\
&=&e^{-x\vee y}\exp\bigl(-W(1)e^{-x\vee y}\bigl(e^x+e^y
\bigr)\bigr)
\\
&&{}\times \exp \biggl(\int_{t=x\vee y}^\infty \P{ \bigl
\{(X_n+x)\wedge (Y_n+y)\ge t \bigr\}}\,
\mathrm{d}t \biggr)\nonumber
\\
&=&e^{-x\vee y}\exp\bigl(-W(1)e^{-x\vee y}\bigl(e^x+e^y
\bigr)\bigr)\nonumber
\\
&&{}\times\exp \bigl(\E{\bigl[\bigl((X_n+x)\wedge(Y_n+y)-x
\vee y\bigr)^+\bigr]} \bigr).\nonumber
\end{eqnarray}
The same calculation also gives
%
%e38 #&#
%e35 ###
\begin{eqnarray}
G(x,y)&=& T^{(2)}\bigl(\mu^{(2)}\bigr) \bigl((x,\infty),(y,
\infty)\bigr)\nonumber
\\
&=&e^{-x\vee y}\exp\bigl(-W(1)e^{-x\vee y}\bigl(e^x+e^y
\bigr)\bigr)
\\
&&{}\times \exp \bigl(\E{\bigl[\bigl((X+x)\wedge(Y+y)-x\vee y\bigr)^+\bigr]} \bigr).\nonumber
\end{eqnarray}

Let
\begin{eqnarray*}
Z_n^{x,y}&:=&\bigl((X_n+x)\wedge(Y_n+y)-x
\vee y\bigr)^+\quad\mbox{and}
\\
Z^{x,y}&:=&\bigl((X+x)\wedge(Y+y)-x\vee y\bigr)^+.
\end{eqnarray*}
Now $(X_n,Y_n)\xrightarrow{\mathrm{D}}(X,Y)$ implies that, for each $(x,y)$,
$Z_n^{x,y}\xrightarrow{\mathrm{D}}Z^{x,y}$. Now
\[
0\le Z_n^{x,y}\le X_n\qquad\mbox{for all }n
\ge1.
\]
Since $\E{X_n}=\E{X}$ for all $n\ge1$, by dominated convergence
theorem, we have $\E{Z_n^{x,y}}\rightarrow \E{Z^{x,y}}$ as
$n\rightarrow\infty$.
Consequently $G_n(x,y)\rightarrow G(x,y)$ for all \mbox{$x,y\in\mathbf{R}_{+}$}.
\end{pf}

%s8.3 #&#
%s8.3 ###
\subsection{Completing the proof of convergence of BP on the PWIT}
With endogeny in hand, we conclude that given a realization of
$\mathcal{T}
$, almost surely, the resulting stationary configuration of the $X$
process of Lemma~\ref{lemtreeProcess} is unique. Also, the following
lemma will show that if the initial messages are i.i.d. random
variables with the fixed point distribution $\mu_*$, then the message
process (\ref{eqBPonT}) converges, and the limit configuration is
unique (almost surely).

%
%le6 #&#
\begin{lemma}\label{lemendL2conv}
If the initial messages $X_{\mathcal{T}}^{0} (\dot{v},v )$ are
i.i.d. random
variables with distribution $\mu_*$, then the message process (\ref
{eqBPonT}) converges in $L^2$ to the process $X$ as $k\rightarrow
\infty$.
\end{lemma}
\begin{pf}
Consider the evolution of bivariate messages according to (\ref
{eqBPonT}), starting from $(X_{\mathcal{T}}^{0} (\cdot
),X  (\cdot ))$. The second
component will remain unchanged because the $X$ process satisfies (\ref
{eqXrecursion}). The distribution of $(X_{\mathcal{T}}^{0}
(\cdot
),X  (\cdot ))$
is $\mu_*\otimes\mu_*$. We have
\[
\operatorname{Law }\bigl(X_{\mathcal{T}}^{k+1} (\cdot ),X  (\cdot
)\bigr)=T^{(2)}\bigl(\operatorname{Law }\bigl(X_{\mathcal{T}}^{k} (
\cdot ),X  (\cdot )\bigr)\bigr).
\]
Here $T^{(2)}$ is as defined in Theorem~\ref{thmendOfECrde}. By
Theorem~\ref{thmendBVU}(c), $(X_{\mathcal{T}}^{k} (\cdot
),X  (\cdot ))$ converges
to $(X  (\cdot ),X  (\cdot
))$ in distribution as $k\rightarrow\infty$. Since
$(X_\mathcal{T}^k-X)^2\le2(X_\mathcal{T}^k)^2+2X^2$, and
$\E{[2(X_\mathcal{T}
^k)^2+2X^2]}=4\E{[X^2]}$, the dominated convergence theorem gives
$\E{[(X_\mathcal{T}^k-X)^2]}\rightarrow0$ as $k\rightarrow\infty$.
\end{pf}

We now prove that if the initial values are i.i.d. random variables
with some arbitrary distribution (not necessarily $\mu_*$), then the
message process (\ref{eqBPonT}) does indeed converge to the unique
stationary configuration. Of course, the initial condition of
particular interest to us is the all zero initial condition (\ref
{eqbpInit}), but we will prove a more general result.

The following lemma will allow us to interchange limit and minimization
while working with the updates on $\mathcal{T}$.

%
%le7 #&#
\begin{lemma} \label{lemcontrolOnT}
Let $X_{\mathcal{T}}^{0} (\dot{v},v )$ be initialized to
i.i.d. random variables
with arbitrary distribution F on $\mathbf{R}_{+}$. Then the map
\[
\pi_{\mathcal{T}}^{k}(v)=\argmin_{u\sim v} \bigl\{ \bigl(\xi
_{\mathcal{T}} (v,u )-X_{\mathcal{T}}^{k} (v,u ) \bigr)^+ \bigr\}
\]
is a.s. well defined and finite for all $k\ge1$, and
\[
\sup_{k\ge1}\P{ \Bigl\{\max\argmin_{i\ge1} \bigl\{
\bigl(\xi _{\mathcal{T}} (v,v.i )-X_{\mathcal{T}}^{k} (v,v.i )
\bigr)^+ \bigr\}\ge i_0 \Bigr\}}\rightarrow0\qquad \mbox{as
}i_0\rightarrow \infty.
\]
\end{lemma}
\begin{pf}
Fix $k$.
If $j\in\argmin_{i\ge1} \{ (\xi_{\mathcal{T}}
(v,v.i )-X_{\mathcal{T}}^{k} (v,v.i ) )^+
\}$ and $j\ge2$, then
\[
\xi (v,v.j )-X_{\mathcal{T}}^{k} (v,v.j )\le \bigl(
\xi_{\mathcal{T}} (v,v.1 )-X_{\mathcal
{T}}^{k} (v,v.1 ) \bigr)^+.
\]
Now
%
%e39 #&#
%e36 ###
\begin{eqnarray}
\label{eqargmin} && \P{ \bigl\{\xi (v,v.j )-X_{\mathcal{T}}^{k}
(v,v.j )\le \bigl(\xi_{\mathcal{T}} (v,v.1 )-X_{\mathcal{T}}^{k}
(v,v.1 ) \bigr)^+ \bigr\}}\nonumber
\\
&&\qquad \le \P{ \bigl\{\xi (v,v.j )\le X_{\mathcal{T}}^{k} (v,v.j ) \bigr
\}}
\\
&&\quad\qquad{}+\P{ \bigl\{\xi (v,v.j )-X_{\mathcal{T}}^{k} (v,v.j )\le\xi
(v,v.1 )-X_{\mathcal{T}}^{k} (v,v.1 ) \bigr\}}.\nonumber
\end{eqnarray}
The\vspace*{1pt} updates are such that $ \{X_{\mathcal{T}}^{k}
(v,v.i ),i\ge1 \}$ remain i.i.d. and independent of the
Poisson process $ \{\xi (v,v.i
) \}$. Thus the
probability on the right-hand side of (\ref{eqargmin}) equals
\[
\P{ \bigl\{\xi_j\le X_1^k \bigr\}}+\P{ \bigl
\{\xi_{j-1}\le X_2^k-X_1^k
\bigr\}},
\]
where $ \{\xi_i \}$ is a Poisson process and $X_1^k,X_2^k$ are
independent random variables with same distribution as $X_{\mathcal{T}
}^{k} (v,v.1 )$.
Then
%
%e40 #&#
%e37 ###
\begin{eqnarray}
\label{eqtailsum} && \sum_{j=2}^\infty \P{
\Bigl\{j\in\argmin_{i\ge1} \bigl\{ \bigl(\xi_{\mathcal{T}} (v,v.i
)-X_{\mathcal
{T}}^{k} (v,v.i ) \bigr)^+ \bigr\} \Bigr\}}\nonumber
\\
&&\qquad \le \sum_{j=2}^\infty \bigl(\P{ \bigl\{
\xi_j\le X_1^k \bigr\} }+\P{ \bigl\{
\xi_{j-1}\le X_2^k-X_1^k
\bigr\}} \bigr)\nonumber
\\
&&\qquad \le \sum_{j=1}^\infty \P{ \bigl\{
\xi_j\le X_1^k \bigr\}}+\sum
_{j=1}^\infty \P{ \bigl\{\xi_{j}\le
X_2^k-X_1^k \bigr\}}
\\
&&\qquad = \E{X_1^k}+\E{\bigl\llvert X_1^k-X_2^k
\bigr\rrvert }\nonumber
\\
&&\qquad \le 3\E{X_1^k}.\nonumber
\end{eqnarray}
From the proof of Theorem~\ref{thmfpConv} it follows that $\E{X_1^k}$ converges, and hence it is bounded. This proves that the
$\argmin$ is a.s. finite and the probability in the statement of the
lemma, being upper bounded by the tail sum of the left-hand side of~(\ref{eqtailsum}), converges uniformly to 0.
\end{pf}

We are now in a position to prove the required convergence.

%
%th13 #&#
\begin{theorem} \label{thmconvOnT}
The recursive tree process defined by (\ref{eqBPonT}) with i.i.d.
initial messages converges to the unique stationary configuration in
the following sense. For every $v\in\mathcal{V}$,
\[
X_{\mathcal{T}}^{k} (v,v.i )\xrightarrow{L^2}X  (v,v.i )\qquad\mbox{as }k\rightarrow\infty.
\]
Also, the decisions at the root converge, that is, $\P{ \{\pi
_{\mathcal{T}}^{k}(\phi)\neq\mathcal{C}_{\mathrm{opt}}(\phi
) \}}\rightarrow0$ as \mbox{$k \rightarrow\infty$}.
\end{theorem}
\begin{pf}
The proof is essentially identical to the proof of Theorem 5.2 of \cite
{SalSha2009}. We present it here for completeness.\vadjust{\goodbreak}

Let $F$ be the c.d.f. of the initial distribution. Let $\theta_t, t\in
\mathbf{R}$ denote the $t$-shift operator on $\mathcal{D}$, that is,
$\theta
_t\widebar{F}\dvtx x\mapsto\widebar{F}(x-t)$. Since $T^n\widebar{F}\rightarrow\widebar{F}_*$,
and $T^n\widebar{F}$ are of the form $y\mapsto c_n e^{-y},y\ge0$ for
$n\ge1$, for any $\varepsilon>0$ there exists $k_\varepsilon\in\mathbf{N}$
such that
\[
\theta_{-\varepsilon}\widebar{F}_*\le T^{k_\varepsilon}\widebar{F}\le\theta
_{\varepsilon}\widebar{F}_*.
\]
By Strassen's theorem, probability measures satisfying such an ordering
can be coupled in a pointwise monotone manner. In other words, there
exists a probability space $E^\prime=(\Omega^\prime,\mathscr
{F}^\prime,P^\prime)$, possibly differing from the original space
$E=(\Omega,\mathscr{F},P)$, on which we can define a random variable\vspace*{2pt}
$X^{\varepsilon}$ with complementary c.d.f. $T^{k_\varepsilon}\widebar{F}$ and
two random variables $X^-$ and $X_+$ with distribution $\widebar{F}_*$
such that almost surely
%
%e41 #&#
%e38 ###
\begin{equation}
\label{eqordering} X^--\varepsilon\le X^\varepsilon\le X^++\varepsilon.
\end{equation}

We now define over the product space $(\bigotimes_{v\in\mathcal{V}
}E^\prime)\otimes E$ the PWIT $\mathcal{T}$ and independent copies
$(X_v^-,X_v^\varepsilon,X_v^+)_{v\in\mathcal{V}}$ of the triple
$(X^-,X^\varepsilon,X^+)$.

On $\mathcal{T}$, we look at the message process with three different
initializations:
\[
X_{\mathcal{T}}^{0,-} (\dot{v},v )=X_v^-,\qquad
X_{\mathcal{T}
}^{0,\varepsilon} (\dot{v},v )=X_v^\varepsilon
\quad\mbox{and}\quad X_{\mathcal{T}}^{0,+} (\dot{v},v )=X_v^+\qquad
\forall v\in\mathcal{V}.
\]
From the update rule (\ref{eqBPonT}) one can readily verify that the
ordering between the messages is preserved in the following sense. For
any $v\in\mathcal{V}$ and $k\ge0$,
\begin{eqnarray*}
X_{\mathcal{T}}^{2k,-} (\dot{v},v )-\varepsilon&\le& X_{\mathcal{T}
}^{2k,\varepsilon}
(\dot{v},v )\le X_{\mathcal{T}}^{2k,+} (\dot{v},v )+\varepsilon;
\\
X_{\mathcal{T}}^{2k+1,+} (\dot{v},v )-\varepsilon&\le& X_{\mathcal{T}
}^{2k+1,\varepsilon}
(\dot{v},v )\le X_{\mathcal{T}
}^{2k+1,-} (\dot{v},v )+\varepsilon.
\end{eqnarray*}

Now fix a $v\in\mathcal{V}$, and observe that
\[
\bigl(X_{\mathcal{T}}^{k+k_\varepsilon} (\dot{v},v )\bigr)_{k\ge
0}\eqdist
\bigl(X_{\mathcal{T}}^{k,\varepsilon} (\dot{v},v )\bigr)_{k\ge0}.
\]
It follows that for every $k\ge k_\varepsilon$,
\begin{eqnarray*}
&& \sup_{s,t\ge k}\bigl\llVert X_{\mathcal{T}}^{s} (
\dot{v},v )-X_{\mathcal{T}}^{t} (\dot{v},v )\bigr\rrVert
_{L^2}
\\
&&\qquad =\sup_{s,t\ge k-k_\varepsilon}\bigl\llVert X_{\mathcal
{T}}^{s,\varepsilon
}
(\dot{v},v )-X_{\mathcal{T}}^{t,\varepsilon} (\dot {v},v )\bigr\rrVert
_{L^2}
\\
&&\qquad \le 2\sup_{t\ge k-k_\varepsilon}\bigl\llVert X_{\mathcal{T}}^{t,\pm}
(\dot {v},v )-X  (\dot{v},v )\bigr\rrVert _{L^2}+2\varepsilon.
\end{eqnarray*}
From endogeny and Lemma~\ref{lemendL2conv}, it follows that
\[
\sup_{t\ge k-k_\varepsilon}\bigl\llVert X_{\mathcal{T}}^{t,\pm} (\dot
{v},v )-X  (\dot{v},v )\bigr\rrVert _{L^2}\rightarrow0\qquad
\mbox{as }k\rightarrow\infty.
\]
Thus the sequence $(X_{\mathcal{T}}^{k} (\dot{v},v
))_{k\ge0}$
is Cauchy in $L^2$,
and hence convergent. Now, Lemma~\ref{lemcontrolOnT} allows us to
interchange limit and minimization in (\ref{eqBPonT}) to conclude
that the limit process has to be a fixed point of (\ref{eqBPonT}). By
endogeny there is a unique stationary configuration a.s. on any
realization of the PWIT. Hence the limit configuration has to be
identical to the $X$ process.

Again by Lemma~\ref{lemcontrolOnT}, for any $\varepsilon>0$, we can
choose an $i_0$ such that
\[
\P{ \bigl\{\pi_{\mathcal{T}}^{k}(\phi)\nsubseteq \{1,2,\ldots,i_0 \} \bigr\}}<\varepsilon/3
\]
for all $k\ge1$, and $\P{ \{\mathcal{C}_{\mathrm{opt}}(\phi
)\nsubseteq \{1,2,\ldots,i_0 \} \}}<\varepsilon/3$.
Now, the convergence of $X_\mathcal{T}^k$
to $X$ implies that for $k$ sufficiently large, when $\pi_{\mathcal
{T}}^{k}(\phi)$ and $\mathcal{C}_{\mathrm{opt}}(\phi)$ are
contained in $ \{1,2,\ldots,i_0 \}$,
the probability that the two maps differ is less than $\varepsilon/3$.
This proves the second statement of the theorem.
\end{pf}

%s9 #&#
%s9 ###
\section{Belief propagation on $\widebar{K}_n$}\label{secbpOnKn}%s9.1 #&#
%s9.1 ###
\subsection{Convergence of the update rule on $\widebar{K}_n$ to the update rule on $\mathcal{T}$}
We use from \cite{SalSha2009} the modified definition of local
convergence applied to geometric networks with edge labels, that is,
networks in which each directed edge $(v,w)$ has a label $\lambda
(v,w)$ taking values in some Polish space. For local convergence of a
sequence of such labeled networks $G_1,G_2,\ldots$ to a labeled
geometric network $G_\infty$, we add the additional requirement that
the rooted graph isomorphisms $\gamma_{n,\rho}$ satisfy
\[
\lim_{n\rightarrow\infty}\lambda_{G_n}\bigl(\gamma_{n,\rho
}(v,w)
\bigr)=\lambda _{G_\infty}(v,w)
\]
for each directed edge $(v,w)$ in $\mathcal{N}_\rho(G_\infty)$.

Now we view the configuration of BP on a graph $G$ at the $k$th
iteration as a labeled geometric network with the label on edge $(v,w)$
given by the pair
\[
\bigl(X_{G}^{k} (v,w ),\mathbf{1}_{ \{v\in\pi
_{G}^{k}(w) \}} \bigr).
\]

With this definition, our convergence result can be written as the
following theorem.

%
%th14 #&#
\begin{theorem} \label{thmupdateConv}
For every fixed $k\ge0$, the $k$th step configuration of BP on
$\widebar{K}_n$
converges in the local weak sense to the $k$th step configuration of BP
on $\mathcal{T}$.
%
%e42 #&#
%e39 ###
\begin{equation}
\label{equpdateConv} \bigl(\widebar{K}_n,X_{\widebar{K}_n}^{k}
(v,w ),\mathbf {1}_{ \{v\in\pi_{\widebar{K}_n}^{k}(w) \}} \bigr)\xrightarrow{\mathrm{l.w.}} \bigl(
\mathcal{T},X_{\mathcal
{T}}^{k} (v,w ),\mathbf{1}_{ \{v\in\pi_{\mathcal
{T}}^{k}(w) \}}
\bigr).
\end{equation}
\end{theorem}
\begin{pf}
The proof of this theorem proceeds along the lines of the proof of
Theorem~4.1 of \cite{SalSha2009}.

Consider an almost sure realization of the convergence $\widebar{K}_n
\rightarrow\mathcal{T}$.

Recall from Section~\ref{seclwConv} the labeling of the vertices of
$\mathcal{T}$ from the set $\mathcal{V}$. We now recursively apply multiple
labels from $\mathcal{V}$ to the vertices of $\widebar{K}_n$. Label
the root as~$\phi$. If $v\in\mathcal{V}$ denotes a vertex $x$ of $\widebar
{K}_n$, then
$(v.1,v.2,\ldots,\break v.(n-1)$ denote the neighbors of $x$ in $\widebar{K}_n$
ordered by increasing lengths of the corresponding edge with~$x$. Then
the convergence in (\ref{equpdateConv}) is shown if we argue that
\begin{eqnarray*}
&&\forall \{v,w \}\in\mathcal{E}\qquad X_{\widebar
{K}_n}^{k} (v,w )
\xrightarrow{\P}X_{\mathcal{T}}^{k} (v,w )\quad\mbox{and}
\\
&&\forall v\in\mathcal{V}\qquad\pi_{\widebar
{K}_n}^{k}(v)
\xrightarrow{\P} \pi_{\mathcal{T}}^{k}(v)\qquad
\mbox{as }n\rightarrow\infty.
\end{eqnarray*}

The above is trivially true for $k=0$. Writing the update and decision
rules as
\begin{eqnarray*}
X_{\widebar{K}_n}^{k+1} (w,v )&=&\min_{u\in \{
v.1,\ldots,v.(n-1),\dot{v} \}\setminus \{w \}
} \bigl\{
\bigl(\xi_{\widebar{K}_n} (v,u )-X_{\widebar
{K}_n}^{k} (v,u ) \bigr)^+
\bigr\}\quad\mbox{and}
\\
\pi_{\widebar{K}_n}^{k}(v)&=&\argmin_{u\in \{v.1,\ldots,v.(n-1),\dot{v} \}} \bigl\{ \bigl(
\xi_{\widebar{K}_n} (v,u )-X_{\widebar{K}_n}^{k} (v,u ) \bigr)^+ \bigr
\},
\end{eqnarray*}
we may try to use the convergence of each term on the right-hand side
inductively to conclude the convergence of the term on the left. This
is not directly possible as the minimum is over an unbounded number of
terms as $n\rightarrow\infty$. However the following lemma allows us to
restrict attention to a uniformly bounded number of terms for each $n$
with probability as high as desired, and hence obtain convergence in
probability for each $k\ge0$.
\end{pf}

%
%le8 #&#
\begin{lemma} \label{lemunifControl}
For all $v\in\mathcal{V}$ and $k\ge0$,
\[
\lim_{i_0\rightarrow\infty}\limsup_{n\rightarrow\infty}\P{ \Bigl\{ \max
\argmin_{1\le i\le n-1} \bigl\{ \bigl(\xi_{\widebar
{K}_n} (v,v.i
)-X_{\widebar{K}_n}^{k} (v,v.i ) \bigr)^+ \bigr\}\ge i_0
\Bigr\}}=0.
\]
\end{lemma}
\begin{pf}
The proof is the same as the proof of Lemma~4.1 of \cite{SalSha2009}.
The only thing to keep in mind is $\argmin$ is a set, and we target
the largest index, but the same proof applies.
\end{pf}
%
%s9.2 #&#
%s9.2 ###
\subsection{Completing the upper bound: Proof of Theorem~\texorpdfstring{\protect\ref{thmbpconv}}{2}}
%We prove Theorem~\ref{thmbpconv} with a slight modification.

%Fix an $M>0$.
By Theorem~\ref{thmconvOnT}, $\pi_{\mathcal{T}}^{k}(\phi
)\xrightarrow{\P}\mathcal{C}_{\mathrm{opt}}(\phi
)\mbox{ as }k\rightarrow\infty$. It follows that
%
%e43 #&#
%e40 ###
\begin{equation}
\label{eqcostconvk} \sum_{v\in\pi_{\mathcal{T}}^{k}(\phi)}\xi_{\mathcal{T}} (\phi,v
) \xrightarrow{\P}\sum_{v\in\mathcal{C}_{\mathrm{opt}}{(\phi)}}
\xi_{\mathcal{T}} (\phi,v )\qquad\mbox{as }k\rightarrow\infty.
\end{equation}
We now prove convergence in expectation.
Observe that
\[
v\in\pi_{\mathcal{T}}^{k}(\phi)\quad \Longrightarrow\quad\xi_{\mathcal
{T}} (
\phi,v )-X_{\mathcal{T}}^{k} (\phi,v )\le \bigl(\xi_{\mathcal{T}} (
\phi,1 )-X_{\mathcal
{T}}^{k} (\phi,1 ) \bigr)^+\le\xi_{\mathcal{T}}
(\phi,1 ).
\]
By (\ref{eqBPonT}), $X_{\mathcal{T}}^{k} (\phi,v )\le
\xi_{\mathcal{T}} (v,v.1 )$. Thus
%
%e44 #&#
%e41 ###
\begin{equation}
\label{eqdominatek} v\in\pi_{\mathcal{T}}^{k}(\phi) \quad \Longrightarrow\quad
\xi_{\mathcal
{T}} (\phi,v )\le\xi_{\mathcal{T}} (\phi,1 )+\xi_{\mathcal{T}}
(v,v.1 ).
\end{equation}
This implies
\[
\sum_{v\in\pi_{\mathcal{T}}^{k}(\phi)}\xi_{\mathcal{T}} (\phi,v )\le
\xi_{\mathcal{T}} (\phi,1 )+\sum_{i\ge2}
\xi_{\mathcal{T}} (\phi,i )\mathbf{1}_{ \{
\xi_{\mathcal{T}} (\phi,i )\le\xi_{\mathcal{T}}
(\phi,1 )+\xi_{\mathcal{T}} (i,i.1 ) \}}.
\]
It can be verified that the sum on the right-hand side in the above
equation is an integrable random variable. Equation~(\ref
{eqcostconvk}) and the dominated convergence theorem give
%
%e45 #&#
\begin{eqnarray}
\label{eqexpCostLimitk} \lim_{k\rightarrow\infty}\E{ \biggl[\sum
_{v\in\pi_{\mathcal
{T}}^{k}(\phi)}\xi_{\mathcal{T}} (\phi,v ) \biggr]} &=&\E{ \biggl[\sum
_{v\in\mathcal{C}_{\mathrm{opt}}{(\phi)}}\xi _{\mathcal{T}} (\phi,v ) \biggr]}
\nonumber\\[-8pt]\\[-8pt]
&=&2W(1)+W(1)^2,\nonumber
\end{eqnarray}
where the last equality follows from Theorem~\ref{thmcoptCost}.

By Theorem~\ref{thmupdateConv} and Lemma~\ref{lemunifControl},
using the definition of local weak convergence, we have
%
%e46 #&#
%e42 ###
\begin{equation}
\label{eqdominaten}
\sum_{v\in\pi_{\widebar{K}_n}^{k}(\phi)}\xi_{\widebar{K}_n}
(\phi,v) \xrightarrow{\P}\sum_{v\in\pi_{\mathcal{T}}^{k}(\phi
)}
\xi_{\mathcal{T}} (\phi,v )\qquad\mbox{as }n\rightarrow\infty.
\end{equation}
We now apply the arguments that lead to (\ref{eqdominatek}) to the
edge covers $\pi_{\widebar{K}_n}^{k}(\phi)$, and obtain
\[
v\in\pi_{\widebar{K}_n}^{k}(\phi) \quad \Longrightarrow\quad\xi_{\widebar
{K}_n} (
\phi,v )\le\xi_{\widebar{K}_n} (\phi,1 )+\xi_{\widebar{K}_n} (v,v.1 ).
\]
For any two vertices $u,v$ of $\widebar{K}_n$, define $S_n(u,v)=\min_{w\neq
u,v}\xi_{\widebar{K}_n} (u,w )$. Then for a vertex $v$ of
$\widebar{K}_n$, $\xi_{\widebar{K}_n} (\phi,1 )\le
S_n(\phi,v)$ and $\xi_{\widebar{K}_n} (v,v.1 )\le
S_n(v,\phi)$. This gives
\[
v\in\pi_{\widebar{K}_n}^{k}(\phi) \quad \Longrightarrow\quad \xi_{\widebar
{K}_n} (
\phi,v )\le S_n(\phi,v)+S_n(v,\phi).
\]
Consequently,
%
%e47 #&#
%e43 ###
\begin{equation}
\label{eqdominaten1} \sum_{v\in\pi_{\widebar{K}_n}^{k}(\phi)}\xi_{\widebar{K}_n} (\phi,v
)\le\sum_v\xi_{\widebar{K}_n} (\phi,v )
\mathbf{1}_{ \{\xi_{\widebar{K}_n} (\phi,v )\le
S_n(\phi,v)+S_n(v,\phi) \}}.
\end{equation}
Observe that $\xi_{\widebar{K}_n} (\phi,v ),S_n(\phi,v)\mbox{ and }S_n(v,\phi)$
are independent exponential random variables with means $n,n/(n-2)$ and
$n/(n-2)$, respectively. So we can write
\begin{eqnarray*}
&& \E{ \bigl[\xi_{\widebar{K}_n} (\phi,v )\mathbf {1}_{ \{\xi_{\widebar{K}_n} (\phi,v )\le S_n(\phi,v)+S_n(v,\phi) \}}
\bigr]}
\\
&&\qquad =\int_0^\infty\int_0^x
\frac{t}{n}e^{-t/n}\, \mathrm{d}t \biggl(\frac
{n-2}{n}
\biggr)^2xe^{- (({n-2})/{n} )x}\, \mathrm{d}x
\\
&&\qquad =\frac{3n^2-5n}{(n-1)^3}.
\end{eqnarray*}
Summing over all neighbors of $\phi$, we get
%
%e48 #&#
%e44 ###
\begin{equation}
\label{eqdominaten2} \E{ \biggl[\sum_v\xi_{\widebar{K}_n}
(\phi,v )\mathbf {1}_{ \{\xi_{\widebar{K}_n} (\phi,v )\le S_n(\phi,v)+S_n(v,\phi) \}} \biggr]}=\frac{3n^2-5n}{(n-1)^2},
\end{equation}
which converges to 3 as $n\rightarrow\infty$.

Using local weak convergence, we can see that
\begin{eqnarray*}
&& \sum_v\xi_{\widebar{K}_n} (\phi,v )
\mathbf{1}_{ \{
\xi_{\widebar{K}_n} (\phi,v )\le S_n(\phi,v)+S_n(v,\phi
) \}}
\\
&&\qquad \xrightarrow{\P}\xi_{\mathcal{T}} (\phi,1 )+\sum
_{i\ge2}\xi_{\mathcal{T}} (\phi,i )\mathbf
{1}_{ \{
\xi_{\mathcal{T}} (\phi,i )\le\xi_{\mathcal{T}}
(\phi,1 )+\xi_{\mathcal{T}} (i,i.1 ) \}}.
\end{eqnarray*}
It can be verified that the expectation of the random variable on the
right-hand side above equals 3. Using this with (\ref{eqdominaten}),
(\ref{eqdominaten1}) and (\ref{eqdominaten2}), the generalized
dominated convergence theorem yields
%
%e49 #&#
%e45 ###
\begin{equation}
\label{eqexpCostLimitn} \lim_{n\rightarrow\infty}\E{ \biggl[\sum
_{v\in\pi_{\widebar{K}_n}^{k}(\phi)}\xi_{\widebar{K}_n} (\phi,v ) \biggr]} =\E{ \biggl[\sum
_{v\in\pi_{\mathcal{T}}^{k}(\phi)}\xi_{\mathcal
{T}} (\phi,v ) \biggr]}.
\end{equation}

Combining (\ref{eqexpCostLimitn}) and (\ref{eqexpCostLimitk}) gives
%
%e50 #&#
%e46 ###
\begin{equation}
\label{eqexpCostLim} \lim_{k\rightarrow\infty}\lim_{n\rightarrow\infty}\E{ \biggl[
\sum_{v\in\pi_{\widebar{K}_n}^{k}(\phi)}\xi_{\widebar{K}_n} (\phi,v ) \biggr]}=
2W(1)+W(1)^2.
\end{equation}

The expectation in the statement of Theorem~\ref{thmbpconv} can be
written as
%
%e51 #&#
%e47 ###
\begin{eqnarray}
\label{eqcostScaling} \E{ \biggl[\sum_{e\in\mathcal{C}(\pi_{K_n}^{k})}
\xi_{K_n} (e ) \biggr]}&=&\frac
{1}{2}\E{ \biggl[\sum
_v\sum_{w\in\pi_{K_n}^{k}(v)}\xi_{K_n }
(v,w ) \biggr]}\nonumber
\\
&=&\frac{1}{2}\E{ \biggl[\sum_v
\frac{1}{n}\sum_{w\in\pi_{\widebar{K}_n}^{k}(v)}\xi_{\widebar{K}_n} (v,w )
\biggr]}
\\
&=&\frac{1}{2}\E{ \biggl[\sum_{v\in\pi_{\widebar{K}_n}^{k}(\phi
)}
\xi_{\widebar{K}_n} (\phi,v ) \biggr]}.\nonumber
\end{eqnarray}
In the first equality above we count the contribution of the edges of
the cover incident at each vertex of $K_n$. The factor of $1/2$
appears because each edge in the edge cover appears twice, once for
each of its endpoints. The $1/n$ in the second equality accounts for
the scaling of edge costs from $K_n$ to $\widebar{K}_n$. The third equality
holds because the root $\phi$ in $\widebar{K}_n$ is chosen uniformly
at random
from the $n$ vertices. Equation~(\ref{eqexpCostLim}) now completes the proof
of Theorem~\ref{thmbpconv}.

%s9.3 #&#
%s9.3 ###
\subsection{Completing the proof of Theorem~\texorpdfstring{\protect\ref{thmlimit}}{1}}
Applying the scaling in (\ref{eqcostScaling}) to the optimal edge
covers in $K_n$ and $\widebar{K}_n$, we get
\[
\E{C_n}=\frac{1}{2}\E{ \biggl[\sum
_{ \{\phi,v \}\in
C_n^*}\xi_{\widebar{K}_n} (\phi,v ) \biggr]}.
\]
Theorem~\ref{thmlowerBd} gives the lower bound
\[
\liminf_{n\rightarrow\infty}\E{C_n}\ge W(1)+\frac{W(1)^2}{2}.
\]
By Theorem~\ref{thmbpconv} for any $\varepsilon>0$, we can find $k$
large such that
\[
\lim_{n\rightarrow\infty}\E{ \biggl[\sum_{e\in\mathcal{C}(\pi
_{K_n}^{k})}
\xi_{K_n} (e ) \biggr]}\le W(1)+\frac{W(1)^2}{2}+\varepsilon.
\]
This gives
\[
\limsup_{n\rightarrow\infty}\E{C_n}\le W(1)+\frac
{W(1)^2}{2}+
\varepsilon.
\]
Since $\varepsilon$ is arbitrary, we get the upper bound
\[
\limsup_{n\rightarrow\infty}\E{C_n}\le W(1)+\frac{W(1)^2}{2}.
\]
This completes the proof of Theorem~\ref{thmlimit}.

Observe that for any $\varepsilon>0$, there exist $K_\varepsilon$ and
$N_\varepsilon$ such that for all $k\ge K_\varepsilon$ and $n\ge N_\varepsilon
$, we have
\[
\E{ \biggl[\sum_{e\in\mathcal{C}(\pi_{K_n}^{k})}\xi_{K_n} (e )
\biggr]}\le W(1)+\frac{W(1)^2}{2}+\varepsilon.
\]
Thus for large $n$ the BP algorithm gives a solution with cost within
$\varepsilon$ of the optimal value in $K_\varepsilon$ iterations. In an
iteration, the algorithm requires $O(n)$ computations at every vertex.
This gives an $O(K_\varepsilon n^2)$ running time for the BP algorithm to
compute an $\varepsilon$-approximate solution. The worst case complexity
of the edge-cover problem is $O(n^3)$, a result due to Edmonds and
Johnson (1970); see \cite{Sch2003}, Theorem~27.2.

%s10 #&#
%s10 ###
\section{More results}\label{secmore}
Our main results for the edge-cover problem were the proof of the limit
of the expected minimum cost (Theorem~\ref{thmlimit}) and the means
to obtain an asymptotically optimal solution using the BP algorithm
(Theorem~\ref{thmbpconv}). The use of objective method as the proof
technique allows us to obtain several auxiliary results about the
structure of the optimal solution, through calculations for the edge
cover $\mathcal{C}_{\mathrm{opt}}$ on the PWIT. In this section we
state and prove, as
examples, results for the distribution of the degree of the root and
the probability that the least cost edge at the root\vadjust{\goodbreak} is part of the
optimal edge cover $\mathcal{C}_{\mathrm{opt}}$. It is easy to show
using local weak
convergence and the results of Sections~\ref{secbp} and \ref{secbpOnKn} that these quantities arise as limits of the quantities
corresponding to the edge covers $\pi_{\widebar{K}_n}^{k}$.

%
%th15 #&#
\begin{theorem}
\[
\P{ \bigl\{\bigl\llvert \mathcal{C}_{\mathrm{opt}}(\phi)\bigr\rrvert =1 \bigr\}
}=e^{-W(1)}\bigl(1+W(1)\bigr).
\]
For $k\ge2$,
\[
\P{ \bigl\{\bigl\llvert \mathcal{C}_{\mathrm{opt}}(\phi)\bigr\rrvert =k \bigr\}
}=e^{-W(1)}\frac{W(1)^k}{k!}.
\]
\end{theorem}
\begin{pf}
As in the proof of Theorem~\ref{thmRDEsolution},
$ \{(\xi_j,X_j),j\ge1 \}$ is a Poisson process on $\mathbf
{R}_{+}\times
\mathbf{R}_{+}$ with intensity $\, \mathrm{d}z \,
\mathrm{d}F_*(x)$.

From the definition of $\mathcal{C}_{\mathrm{opt}}$,
\begin{eqnarray*}
\P{ \bigl\{\bigl\llvert \mathcal{C}_{\mathrm{opt}}(\phi)\bigr\rrvert =1 \bigr\}
}&=&\P{ \bigl\{ \mbox{at most one point of } \bigl\{(\xi_j,X_j)
\bigr\}}
\mbox{ in } \bigl\{(z,x)\dvtx z-x\le0 \bigr\} \bigr\}
\\
&=&e^{-A}(1+A),
\end{eqnarray*}
where
\begin{eqnarray*}
A&=&\int_{z=0}^{\infty}\int_{x=z}^{\infty}
\, \mathrm{d}F_*(x)\, \mathrm{d}z
\\
&=&\int_{z=0}^{\infty}W(1)e^{-z}\,
\mathrm{d}z
\\
&=&W(1).
\end{eqnarray*}
Thus
\[
\P{ \bigl\{\bigl\llvert \mathcal{C}_{\mathrm{opt}}(\phi)\bigr\rrvert =1 \bigr\}
}=e^{-W(1)}\bigl(1+W(1)\bigr).
\]

For $k\ge2$,
\begin{eqnarray*}
\P{ \bigl\{\bigl\llvert \mathcal{C}_{\mathrm{opt}}(\phi )\bigr
\rrvert =k \bigr\}}&=&\P{ \bigl\{ k\mbox{ points of } \bigl\{(\xi_j,X_j)
\bigr\}\mbox{ in } \bigl\{ (z,x)\dvtx z-x\le0 \bigr\} \bigr\}}
\\
&=&e^{-A}\frac{A^k}{k!}
\\
&=&e^{-W(1)}\frac{W(1)^k}{k!}.
\end{eqnarray*}
\upqed
\end{pf}

%
%th16 #&#
\begin{theorem}
\[
\P{ \bigl\{1\in\mathcal{C}_{\mathrm{opt}}(\phi) \bigr\}}=\frac
{W(1)}{2}+
\frac{1}{W(1)}-W(1)^2-1.
\]
\end{theorem}
\begin{pf}
The event $ \{1\in\mathcal{C}_{\mathrm{opt}}(\phi) \}$
equals the union of two disjoint events:
\begin{longlist}[(a)]
\item[(a)] $\xi (\phi,1 )-X  (\phi,1
)<0$ and
\item[(b)] $0\le\xi (\phi,1 )-X  (\phi,1 )\le\xi (\phi,i )-X  (\phi,i )$
for all $i\ge2$.\vadjust{\goodbreak}
\end{longlist}

The probability of the first event is
\begin{eqnarray*}
\P{ \bigl\{\xi (\phi,1 )-X  (\phi,1 )<0 \bigr\}}&=&\int
_{z=0}^\infty\int_{x=z}^\infty
\, \mathrm{d}F_*(x)e^{-z}\, \mathrm{d}z
\\
&=&\int_{z=0}^\infty W(1)e^{-z}e^{-z}
\, \mathrm{d}z
\\
&=&\frac{W(1)}{2}.
\end{eqnarray*}

For the second event, write $\xi (\phi,i )=\xi (\phi,1 )+\xi
_i^\prime$, where obviously $ \{\xi_i^\prime,i\ge2 \}$
is a rate 1
Poisson process independent of $ \{X  (\phi,i ),i\ge2 \}$. For $i\ge
2$, $\xi (\phi,1 )-X  (\phi,1 )\le
\xi (\phi,i )-X  (\phi,i )$ if and
only if $-X  (\phi,1 )\le\xi_i^\prime
-X  (\phi,i )$. The probability
of the second event can be written as
\begin{eqnarray*}
&& \P{ \bigl\{0\le\xi (\phi,1 )-X  (\phi,1 )\le\xi
(\phi,i )-X  (\phi,i )\mbox{ for all }i\ge2 \bigr\} }
\\[-1pt]
&&\qquad =\int_{x_1=0}^\infty\int_{z_1=x_1}^\infty
P \bigl\{\mbox{no point of } \bigl\{\bigl(\xi_i^\prime,X
  (\phi,i ),i\ge2\bigr) \bigr\}
\\[-1pt]
&&\hspace*{135pt} \mbox{ in } \bigl
\{(z,x)\dvtx z-x\le-x_1 \bigr\} \bigr\} e^{-z_1}\,
\mathrm{d}z_1\, \mathrm{d} F_*(x_1)
\\[-1pt]
&&\qquad =\int_{x_1=0}^\infty e^{-x_1} \exp \biggl(-
\int_{z=0}^\infty\int_{x=z+x_1}^\infty
\, \mathrm{d}F_*(x)\, \mathrm{d}z \biggr) \,
\mathrm{d}F_*(x_1)
\\[-1pt]
&&\qquad =\int_{x_1=0}^\infty e^{-x_1} \exp \biggl(-
\int_{z=0}^\infty W(1)e^{-z}e^{-x_1}
\, \mathrm{d}z \biggr) \, \mathrm{d}F_*(x_1)
\\[-1pt]
&&\qquad =\int_{x_1=0}^\infty e^{-x_1} \exp
\bigl(-W(1)e^{-x_1} \bigr) \, \mathrm{d} F_*(x_1)
\\[-1pt]
&&\qquad =W(1) \bigl(1-W(1)\bigr)+\int_{x_1=0}^\infty W(1)
e^{-2x_1} \exp \bigl(-W(1)e^{-x_1} \bigr) \,
\mathrm{d}x_1
\\[-1pt]
&&\qquad =W(1) \bigl(1-W(1)\bigr)+\frac{1}{W(1)}-W(1)-1
\\[-1pt]
&&\qquad =\frac{1}{W(1)}-W(1)^2-1.
\end{eqnarray*}\upqed
\end{pf}
%
%s11 #&#
%s11 ###
\section{Summary} \label{secsummary}
In a nutshell, we have implemented Aldous's program based on \cite
{Ald2001} to solve the random edge-cover problem. Aldous's program
serves as a rigorous mathematical alternative to the cavity method
applied to mean-field combinatorial optimization problems. Aldous and
Bandyopadhyay \cite{AldBan2005}, Section~7.5, outline the steps of this
rigorous methodology, highlighting the role of RDEs and endogeny. See below.

But first, we must indicate another way in which the complete graph
with i.i.d. edge weights arises. Combinatorial optimization problems
involving $n$ random points on $\mathbb{R}^d$ are of interest in many
physical settings, but are typically difficult to analyze because of
dependence of the random variables representing the $\binom{n}{2}$
distances. A more tractable \textit{mean-field model} ignores the
underlying $d$-dimensional space, and simply models the interpoint
distances as i.i.d. random variables. This resulting model is then the
complete graph on $n$ vertices with i.i.d. edge weights. The case of
exponential mean 1 edge weights models the $d=1$ setting. There are
other distributions to model the $d > 1$ settings. Though we did not
deal with $d > 1$ in this paper, we expect the extension to hold (as
for matching).

Let us return to Aldous's program, as summarized by Aldous and
Bandyopadhyay \cite{AldBan2005}, Section~7.5, and reproduced below.\vspace*{6pt}
\begin{eqnarray*}
&& \mbox{``Start with a combinatorial optimization problem over some}
\\[-6pt]
&& \mbox{size-$n$ random structure.}
\end{eqnarray*}

\begin{itemize}
\item Formulate a ``size-$\infty$'' random structure, the $n
\rightarrow
\infty$ limit in the sense of local weak convergence.
\item Formulate a corresponding combinatorial optimization problem on
the size-$\infty$ structure.
\item Heuristically define relevant quantities on the size-$\infty$
structure via additive renormalization \ldots
\item If the size-$\infty$ structure is treelike (the only case where
one expects exact asymptotic solutions), observe that the relevant
quantities satisfy a problem dependent RDE.
\item Solve the RDE. Use the unique solution to find the value of the
optimization problem on the size-$\infty$ structure.
\item Show that the RTP associated with the solution is endogenous.
\item Endogeny shows that the optimal solution is a measurable function
of the data, in the infinite-size problem. Since a measurable function
is almost continuous, we can pull back to define almost-feasible
solutions of the size-$n$ problem with almost the same cost.
\item Show that in the size-$n$ problem one can patch an
almost-feasible solution into a feasible solution for asymptotically
negligible cost.'' \cite{AldBan2005}, Section~7.5.
\end{itemize}

The size-$n$ random structure is the complete graph on $n$-vertices
$\widebar{K}_n$ with independent exponential mean-$n$ edge weights.
The following points elaborate on how we addressed the steps above:

\begin{itemize}
\item The size-$\infty$ random structure is the PWIT.
\item The corresponding optimization problem on the size-$\infty$
structure is simply the minimum-cost edge cover on the PWIT. While this
step is easy for the edge-cover problem, in general some subtleties are
involved. For example, the limiting size-$\infty$ problem for Frieze's
size-$n$ problem of minimal spanning tree on $\widebar{K}_n$ \cite
{Fri1985} is a minimal spanning forest with certain requirements on the
included edges. See \cite{AldSte2004}, Definition~4.2, for details.
\item We then heuristically provided the quantities relevant to the
edge-cover problem on the PWIT in Section~\ref{secrde}. The additive
renormalization measured the reduction in cost arising from the
relaxation of the requirement that the root be hit.
\item Using the tree structure of the limiting object, we obtained the
RDE (\ref{eqECrde}) associated with the edge-cover problem.
\item We solved the RDE in Theorem~\ref{thmRDEsolution}, showed that
it had a unique solution, and found the value of the optimization
problem on the PWIT in Theorem~\ref{thmcoptCost}. Another important
step is Theorem~\ref{thmCopt} which proves that the edge cover
$\mathcal{C}_{\mathrm{opt}}
$, based on the heuristic relation (\ref{eqDT}), is optimal among
involution invariant edge covers on the PWIT. Our method for
establishing this nontrivial step may have some bearing on other
similar combinatorial optimization problems. This step eventually
established a lower bound for the liminf of size-$n$ optimal values.
\item Theorem~\ref{thmendOfECrde} established endogeny of the RTP
associated with the solution of~(\ref{eqECrde}). Theorem~\ref
{thmbpconv} corresponding to the BP algorithm on $K_n$ replaces the
procedure of Aldous's program for obtaining solutions of the size-$n$
problem from the solution of the size-$\infty$ problem. The key steps
for this are based on Salez and Shah's approach \cite{SalSha2009} and
is as follows. Using endogeny, we argued that BP (with i.i.d.
initializations) converges to the RDE-based stationary configuration on
the PWIT. We then established that, at a particular node of $\widebar{K}_n$,
the BP update for large $n$ depends essentially only on messages from
its local neighborhood (Lemma~\ref{lemunifControl}). This is then
used to express BP on the PWIT as the limit of BP on $\widebar{K}_n$.
The BP
iterates on $\widebar{K}_n$ were then the candidate solutions for the
size-$n$ problem.
\item No corrective patch-up was needed for the size-$n$ problem, since
at each iteration of the BP algorithm, every vertex was covered by the
corresponding selection of edges. Simple dominated convergence
arguments then established the convergence of the expected optimal
costs to the correct value.
\end{itemize}

It is worth noting that the upper bound result in Theorem~\ref
{thmlimit} can be obtained via a simpler proof of Theorem~\ref
{thmbpconv} for a version of BP algorithm, where the messages are
initialized as i.i.d. random variables from the fixed-point
distribution $F_*$. In this case Lemma~\ref{lemendL2conv}, which
follows from endogeny, establishes the convergence result on the PWIT.
The more general result of Theorem~\ref{thmconvOnT} shows that BP
works when messages are initialized as i.i.d. random variables from
any arbitrary distribution.

Finally, we must mention that Aldous \cite{Ald2001} proved a strong
property called \textit{asymptotic essential uniqueness} for matching,
which is roughly the property that if a matching on $\widebar{K}_n$
is almost optimal, then it coincides with the optimal matching, except
on a small proportion of edges. The question of whether this property
holds for the edge-cover problem is one that we hope to address in the
near future.

\section*{Acknowledgments}
Part of this work was carried out when
Rajesh Sundaresan was on sabbatical leave at the University of
Illinois at Urbana--Champaign whose support is gratefully
acknowledged.

%suskaldyti doi

% imsref loaded by linak, 2014-03-14 18:14:41
%
% imsref loaded by linak, 2014-03-17 15:05:37

\printaddresses

\end{document}